\documentclass[11pt]{amsproc}

\usepackage{eufrak}
\usepackage{amscd}

\theoremstyle{plain}
\newtheorem{thm}{Theorem}[section]
\newtheorem{lem}[thm]{Lemma}
\newtheorem{prop}[thm]{Proposition}
\newtheorem{cor}[thm]{Corollary}
\newtheorem*{ta}{Theorem 1}
\newtheorem*{tb}{Theorem 2}
\newtheorem*{tc}{Theorem 3}

\newcommand{\pl}{$\Gamma = ({\mathcal P}, {\mathcal L})$}
\newcommand{\nas}[1]{N_{\mathcal A}(#1)}
\newcommand{\mss}[2]{M_{#1}{(#2)}}
\newcommand{\res}[2]{Res_{#1}{(#2)}}

\begin{document}

\title[]{A characterization of two classes of locally truncated diagram geometries}
\maketitle

\begin{center}
\bigskip
{\Large Silvia Onofrei\footnote{Email address: onofrei@math.ksu.edu.}}\\
{\it Department of Mathematics, Kansas State University,\\ 138 Cardwell
Hall, Manhattan, KS 66506}

\medskip
{\it Advances in Geometry, (4)2004, 469-495. }\\
\end{center}

\begin{abstract} Let {\pl} be a parapolar space which is locally
$A_{n-1,3}({\mathbb K})$ for some integer $n > 6$ and ${\mathbb K}$ a field. There exists a class of
$2$-convex subspaces ${\mathbb D}$, each isomorphic to $D_{5,5}({\mathbb K})$, such
that every symplecton of $\Gamma$ is contained in a unique element of ${\mathbb
D}$. Let {\pl} be a parapolar space which is
locally $A_{n-1,4}({\mathbb K})$ for $n = 7$ or an integer $n \geq 9$ and some field ${\mathbb K}$. Assume that $\Gamma$ satisfies the extra condition called the Weak Hexagon Axiom. Then there exists a class of $2$-convex subspaces ${\mathcal D}$, each isomorphic to
$D_{6,6}({\mathbb K})$, such that every symplecton of $\Gamma$ is contained in a unique
element of ${\mathcal D}$. In both of the above cases $\Gamma$ is the homomorphic image of a truncated building.
\end{abstract}


\bigskip
\section{Introduction}

The purpose of this paper is to give a characterization of two families of locally truncated diagram geometries which arise in the study of the geometries related to the exceptional buildings $E_n$ of both spherical and affine type.\\

Let $I = \lbrace 1, \ldots n \rbrace$ be a finite index set and $K'$ (respectively $K$) a proper subset of $I$. Set $J' = I \setminus K'$ and $J = I \setminus K$. Assume first that {\pl} is a parapolar space which is locally $A_{n-1, 3}({\mathbb K})$ where $n$ is an integer greater than $6$ and ${\mathbb K}$ is a field. There are two classes of maximal singular subspaces denoted by ${\mathcal A}$ and ${\mathcal B}$. Cooperstein theory provides us with a class of symplecta ${\mathcal S}$, which are nondegenerate polar spaces of type $D_4$. Therefore $\Gamma$ is a rank four geometry $\Gamma  =({\mathcal B}, {\mathcal P}, {\mathcal L}, {\mathcal A})$ over $K' =\lbrace 2, \ldots 5 \rbrace$ which is $J'$-locally truncated with diagram:\\
{\it 
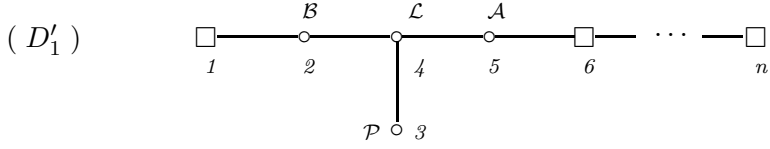
\begin{picture}(1000, 50)(0,0)
\put(0,25){$( \; D'_1 \; )$}
\put(60,25){$\qed$}
\put(75,15){\scriptsize 1}
\put(80,29){\line(1,0){30}}
\put(110,26){$\circ$}
\put(112,36){\scriptsize ${\mathcal B}$}
\put(112,15){\scriptsize 2}
\put(115,29){\line(1,0){30}}
\put(145,26){$\circ$}
\put(152,36){\scriptsize ${\mathcal L}$}
\put(154,15){\scriptsize 4}
\put(148,-3){\line(0,1){30}}
\put(145,-9){$\circ$}
\put(135,-10){\scriptsize ${\mathcal P}$}
\put(154,-10){\scriptsize 3}
\put(150,29){\line(1,0){30}}
\put(180,26){$\circ$}
\put(182,36){\scriptsize ${\mathcal A}$}
\put(182,15){\scriptsize 5}
\put(185,29){\line(1,0){30}}
\put(203,25){$\qed$}
\put(218,15){\scriptsize 6}
\put(223,29){\line(1,0){15} \; $\dots$  \line(1,0){15}}
\put(268,25){$\qed$}
\put(283, 15){\scriptsize n}
\end{picture}}
\vspace{.3cm}

We prove the following:
\begin{ta} Let {\pl} be a parapolar space which is locally $A_{n-1,3}({\mathbb K})$
for some integer $n > 6$ and field ${\mathbb K}$. There is a
collection ${\mathbb D}$ of $2$-convex subspaces, each isomorphic to
$D_{5,5}({\mathbb K})$ such that, if $S \in {\mathcal S}$ is a symplecton, there exists a unique member ${\mathsf D}(S) \in {\mathbb D}$ containing $S$.
\end{ta}
Theorem $1$ says that the above parapolar space {\pl} is enriched to a rank five geometry $\Gamma = ({\mathcal B}, {\mathcal P}, {\mathcal L}, {\mathcal A}, {\mathbb D})$ over $K = \lbrace 2, \ldots 6 \rbrace$ which is $J$-locally truncated with diagram:\\
{\it 
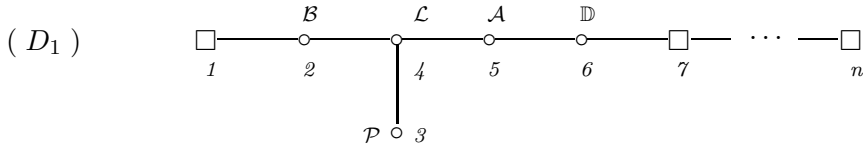
\begin{picture}(1000, 50)(0,0)
\put(0,25){$( \; D_1 \; )$}
\put(60,25){$\qed$}
\put(75,15){\scriptsize 1}
\put(80,29){\line(1,0){30}}
\put(110,26){$\circ$}
\put(112,36){\scriptsize ${\mathcal B}$}
\put(112,15){\scriptsize 2}
\put(115,29){\line(1,0){30}}
\put(145,26){$\circ$}
\put(154,36){\scriptsize ${\mathcal L}$}
\put(154,15){\scriptsize 4}
\put(148,-3){\line(0,1){30}}
\put(145,-9){$\circ$}
\put(135,-10){\scriptsize ${\mathcal P}$}
\put(154,-10){\scriptsize 3}
\put(150,29){\line(1,0){30}}
\put(180,26){$\circ$}
\put(182,36){\scriptsize ${\mathcal A}$}
\put(182,15){\scriptsize 5}
\put(185,29){\line(1,0){30}}
\put(215,26){$\circ$}
\put(217,36){\scriptsize ${\mathbb D}$}
\put(217,15){\scriptsize 6}
\put(220,29){\line(1,0){30}}
\put(239,25){$\qed$}
\put(253,15){\scriptsize 7}
\put(259,29){\line(1,0){15} \; $\dots$  \line(1,0){15}}
\put(304,25){$\qed$}
\put(319,15){\scriptsize n}
\end{picture}}
\vspace{.3cm}

Next, let {\pl} be a parapolar space which is locally $A_{n-1, 4}({\mathbb K})$ where $n$ is an integer, $n=7$ or $n \geq 9$ and ${\mathbb K}$ is a field. Furthermore, assume that $\Gamma$ satisfies the following:\\

{\bf The Weak Hexagon Axiom} ({\it WHA}): {\it Let $H = (x_1, \ldots x_6)$ be a
$6$-circuit, isometrically embedded in $\Gamma$, this means that $x_i \in x^{\perp}_{i+1}$, indices taken mod $6$, and all the other pairs are not collinear. Also assume that at least one of the pairs of points at distance two, say $\lbrace x_1, x_3 \rbrace$, is polar. Then there exists a point $w \in x_1^{\perp} \cap x_3^{\perp} \cap x_5^{\perp}$}.\\

From the local properties it follows that there are two classes of maximal singular subspaces ${\mathcal A}$ and ${\mathcal B}$ and a class ${\mathcal S}$ of symplecta, which are nondegenerate polar spaces of type $D_4$. Therefore $\Gamma = ({\mathcal B}, {\mathcal P}, {\mathcal L}, {\mathcal A})$ is a rank four geometry over $K' = \lbrace 3, \ldots 6 \rbrace$ which is $J'$-locally truncated with the following diagram:\\
{\it 
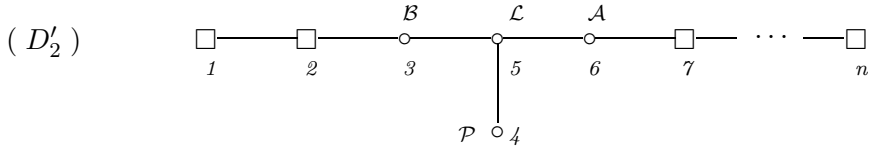
\begin{picture}(1000, 50)(0,0)
\put(0,25){$( \; D'_2 \; )$}
\put(60,25){$\qed$}
\put(75,15){\scriptsize 1}
\put(80,29){\line(1,0){30}}
\put(98,25){$\qed$}
\put(113,15){\scriptsize 2}
\put(118,29){\line(1,0){30}}
\put(148,26){$\circ$}
\put(150,36){\scriptsize ${\mathcal B}$}
\put(150,15){\scriptsize 3}
\put(153,29){\line(1,0){30}}
\put(183,26){$\circ$}
\put(190,36){\scriptsize ${\mathcal L}$}
\put(190,15){\scriptsize 5}
\put(186,-3){\line(0,1){30}}
\put(183,-9){$\circ$}
\put(171,-10){\scriptsize ${\mathcal P}$}
\put(190,-10){\scriptsize 4}
\put(188,29){\line(1,0){30}}
\put(218,26){$\circ$}
\put(220,36){\scriptsize ${\mathcal A}$}
\put(220,15){\scriptsize 6}
\put(223,29){\line(1,0){30}}
\put(241,25){$\qed$}
\put(255,15){\scriptsize 7}
\put(261,29){\line(1,0){15} \; $\dots$  \line(1,0){15}}
\put(306,25){$\qed$}
\put(321, 15){\scriptsize n}
\end{picture}}
\vspace{.3cm}

In Section $5$ we prove:
\begin{tb} Let {\pl} be a parapolar space which is locally $A_{n-1,4}({\mathbb K})$
for some integer $n=7$ or $n \geq 9$ and field ${\mathbb K}$. Assume $\Gamma$ satisfies (WHA). Then there is a collection ${\mathcal D}$ of $2$-convex subspaces, each isomorphic to the
half-spin geometry $D_{6,6}({\mathbb K})$, such that, if $S \in {\mathcal S}$ is a
symplecton, there exists a unique member $D(S) \in {\mathcal D}$
containing $S$.
\end{tb}
The Theorem $2$ provides us with a new class of objects ${\mathcal D}$ and, consequently, $\Gamma = \lbrace {\mathcal B}, {\mathcal P}, {\mathcal L}, {\mathcal A}, {\mathcal D} \rbrace$ is a rank five geometry over $K = \lbrace 3, \ldots 7 \rbrace$ which is $J$-locally truncated with the following diagram:\\
{\it 
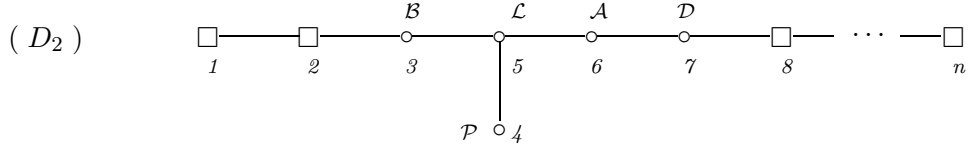
\begin{picture}(1000, 50)(0,0)
\put(0,25){$( \; D_2 \; )$}
\put(60,25){$\qed$}
\put(75,15){\scriptsize 1}
\put(80,29){\line(1,0){30}}
\put(98,25){$\qed$}
\put(113,15){\scriptsize 2}
\put(118,29){\line(1,0){30}}
\put(148,26){$\circ$}
\put(150,36){\scriptsize ${\mathcal B}$}
\put(150,15){\scriptsize 3}
\put(153,29){\line(1,0){30}}
\put(183,26){$\circ$}
\put(190,36){\scriptsize ${\mathcal L}$}
\put(190,15){\scriptsize 5}
\put(186,-3){\line(0,1){30}}
\put(183,-9){$\circ$}
\put(171,-10){\scriptsize ${\mathcal P}$}
\put(190,-10){\scriptsize 4}
\put(188,29){\line(1,0){30}}
\put(218,26){$\circ$}
\put(220,36){\scriptsize ${\mathcal A}$}
\put(220,15){\scriptsize 6}
\put(223,29){\line(1,0){30}}
\put(253,26){$\circ$}
\put(253,36){\scriptsize ${\mathcal D}$}
\put(255,15){\scriptsize 7}
\put(258,29){\line(1,0){30}}
\put(277,25){$\qed$}
\put(292,15){\scriptsize 8}
\put(297,29){\line(1,0){15} \; $\dots $ \line(1,0){15}}
\put(342,25){$\qed$}
\put(357,15){\scriptsize n}
\end{picture}}
\vspace{.3cm}

In Section $6$, using Brouwer-Cohen sheaf theory \cite{loc} and a result of Ellard and Shult \cite{els}, we construct a residually connected sheaf over a locally truncated
geometry with the diagram:\\
{\it 
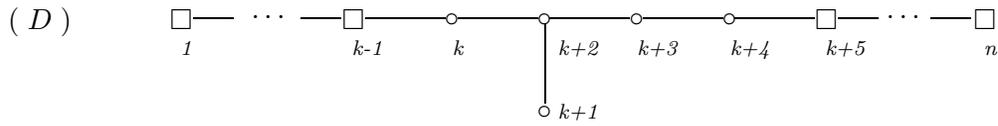
\begin{picture}(1000, 50)(0,0)
\put(0, 25){$(\; D\;)$}
\put(50,25){$\qed$}
\put(70,29){\line(1,0){15}  \; $\dots$ \line(1,0){15}}
\put(65,15){\scriptsize 1}
\put(115,25){$\qed$}
\put(130,15){\scriptsize k-1}
\put(135,29){\line(1,0){30}}
\put(165,26){$\circ$}
\put(168,15){\scriptsize k}
\put(170,29){\line(1,0){30}}
\put(200,26){$\circ$}
\put(208,15){\scriptsize k+2}
\put(203,-3){\line(0,1){30}}
\put(200,-9){$\circ$}
\put(208,-9){\scriptsize k+1}
\put(205,29){\line(1,0){30}}
\put(235,26){$\circ$}
\put(238,15){\scriptsize k+3}
\put(240,29){\line(1,0){30}}
\put(270,26){$\circ$}
\put(273,15){\scriptsize k+4}
\put(275,29){\line(1,0){30}}
\put(294,25){$\qed$}
\put(309,15){\scriptsize k+5}
\put(314,29){\line(1,0){15} \;$\ldots$ \line(1,0){15}}
\put(354,25){$\qed$}
\put(369,15){\scriptsize n}
\end{picture}}
\vspace{.3cm}

The existence of a sheaf for such a locally
truncated geometry was stated without proof in Brouwer and Cohen (Theorem $4$, \cite{loc}). For completeness, since the author does not know of any place in the literature where this sheaf was constructed, we shall give a detailed description of this sheaf.\\

In Section $7$, combining the sheaf theoretic approach with the functorial relation between geometries and chamber systems and with Tits' Local Approach Theorem, we get the following:
\begin{tc} {\bf a.} Let {\pl} be a parapolar space which is locally $A_{n-1,3}({\mathbb K})$ where $n$ is an integer greater than $6$ and ${\mathbb K}$ is a field. Then $\Gamma$ is a residually connected locally truncated diagram geometry belonging to the diagram $(D_1)$ whose universal $2$-cover is the truncation of a building.\\
{\bf b.} Let {\pl} be a parapolar space which is locally $A_{n-1, 4}({\mathbb K})$ where $n$ is an integer $n=7$ or $n \geq 9$ and ${\mathbb K}$ is a field. Assume that $\Gamma$ satisfies the Weak Hexagon Axiom. Then $\Gamma$ is a residually connected diagram geometry belonging to the diagram $(D_2)$ whose universal $2$-cover is the truncation of a building.
\end{tc}

{\it Remarks}:\\
$\bf {1.}$ The diagrams $(D'_1)$ and $(D_1)$ are also denoted $E_n,\; n\geq 4$, see Pasini (Exercise $5.36$ in \cite{pas}).\\

$\bf {2.}$ In Theorem $1$, $n$ is assumed to be at least $7$. When $n = 5$, $\Gamma$ itself is a half-spin geometry of type $D_{5,5}$. When $n=6$ all the maximal singular subspaces at a point have the same rank and we cannot recognize a global partition in two classes according to dimension. In this case $\Gamma \simeq E_{6,4}$ which was characterized by Cohen and Cooperstein \cite{cc}.\\

$\bf {3.}$ The cases $n=7$ and $n=8$ for the diagram $(D_1)$ give the Coxeter diagrams of type $E_7$ and $E_8$. Their geometries, of exceptional type, have also been studied by Hanssens \cite{hans}.\\

$\bf {4.}$ In the hypotheses of Theorem $2$ we assume $n=7$ or $n \geq 9$. When $n=6$, $\Gamma$ is a half-spin geometry of type $D_{6,6}$. The case $n=8$ is excluded for the reason that all the maximal singular subspaces at a point have the same rank and therefore a global partition in two classes is not obvious anymore. However, one can still prove the existence of the $2$-convex subspaces of type $D_{6,6}$. Then every symplecton lies in at least two such subspaces; the details of the proof can be found in Onofrei (Section $3.4$ in \cite{thesis}). This geometry is related to the affine building of type $\widetilde E _7$.\\

$\bf{5.}$ The Weak Hexagon Axiom was used, by our knowledge, at least twice before in the literature. In Shult \cite{shm}, a version of ({\it WHA}) is mentioned under the name of ``the hexagon property" in connection with geometries related to extended $F_4$. Also El-Atrash and Shult \cite{als} used a ``weak hexagon property" for the case of strong parapolar spaces. Their ``weak hexagon property" is the same as ({\it WHA}) but adapted to the case of strong parapolar spaces.\\

$\bf {6.}$ We have assumed $n$ to be finite. However, the results of the Theorems $1$ and $2$ can be extended to geometries of arbitrary rank, since the proofs involve only the class ${\mathcal A}$ of maximal singular subspaces, which can have singular rank at most $5$, and truncations of small rank of the other class ${\mathcal B}$. In Theorem $3$, the finiteness of $n$ is essential since otherwise the sheaf and its associated chamber system are not residually connected anymore; Kasikova and Shult \cite{ks2}.

\section{Preliminaries and definitions}

We assume the reader is familiar with the basic definitions related to point-line geometries. A standard reference is {\it The Handbook of Incidence Geometry} \cite{hnbk}, Chapters $3$ and $12$.

\subsection{Geometries} Let $\Gamma$ be a geometry over $I$, that is a system $\Gamma = (V, *, t)$ consisting of set $V$, a binary, symmetric, reflexive relation on $V$ and a mapping $t: V \rightarrow I$. The elements of $V$ are called {\it objects}, $*$ is called {\it incidence relation} and $t$ is the {\it type function} of $\Gamma$.

\medskip
A {\it flag} $F$ of $\Gamma$ is a (possibly empty) subset of pairwise incident objects of $\Gamma$. The set $t(F)$ is the {\it type} of $F$ and the set $I \setminus t(F)$ is its {\it cotype}. The cardinalities of these sets are the {\it rank} and the {\it corank} of $F$. The {\it residue} of $F$ in $\Gamma$ is the geometry $\res{\Gamma}{F} = (V_F, *_{|V_F}, t _{|V_F})$ over $I \setminus t(F)$, where $V_F$ is the set of all members of $V \setminus F$ incident with each element of $F$. The corank of a flag $F$ is the rank of $\res {\Gamma}{F}$.

\medskip
Let $\Gamma _k = (V_k, *_k, t_k)$ with $k \in K$ some index set and where each $\Gamma_k$ is a geometry over $I_k$. Assume that $\lbrace I_k \rbrace _{k \in K}$ is a family of pairwise disjoint sets. We define the {\it direct sum of geometries}, the geometry denoted by $\Gamma = \bigoplus _{k \in K} \Gamma _k = (V, *, t)$, where $V = \bigcup _{k \in K} V_k$. The incidence is defined as follows: $*_{|V_k}\;  :=  *_k$ and $x * y$ for any two objects $x \in V_{k_1}$ and $y \in V_{k_2}$ with $k_1 \not = k_2$. Finally $t _{|V_k}\;  := t_k$. \\

A geometry $\Gamma$ is {\it residually connected} if and only if for every flag $F$ of corank at least one, $\res{\Gamma}{F}$ is not empty and if, for each flag of corank at least two, $\res{\Gamma}{F}$ is connected.\\

We can define a category ${\mathcal Geom}_I$, in which the objects are
geometries with typeset  $I$ with their corresponding morphisms, i.e. the
type preserving graph morphisms.\\

\subsection{Locally truncated geometries}
For a more detailed account of the concepts from this subsection the reader is referred to Brouwer and Cohen \cite{loc} and Ronan \cite{ron}. Let $I$ be an index set and let $J \subset I$. The {\it truncation of type $J$} of $\Gamma$, denoted
by $^J \Gamma$ is the geometry obtained by restricting the typeset of
$\Gamma$ to $J$. The truncation is a functor from the category of geometries over $I$
to the category of geometries over $J$:
\begin{center}
$^J {\bf Tr}: {\mathcal Geom}_I \rightarrow {\mathcal Geom}_J$.
\end{center}
The {\it $J$-truncation of $\Gamma$}, denoted by $_J \Gamma$, has as objects the objects of $\Gamma$ whose types
are in $I \setminus J$, incidence and type function are those from
$\Gamma$ but restricted to $I \setminus J$. Differently said, the
$J$-truncation of $\Gamma$ is the truncation of type $I \setminus J$ of $\Gamma$: $_J \Gamma = \; ^{I \setminus J}\Gamma$.\\

A {\it diagram} $D$ over $I$ is a mapping which assigns to each $2$-subset $\lbrace i, j \rbrace$ of $I$, a class $D(i,j)$ of rank $2$ geometries.  A geometry $\Delta$ over $I$ {\it belongs to the diagram} $D$ if and only if every
residue of type $\lbrace i, j \rbrace$ of $\Delta$ is a geometry from $D(i,j)$.\\

A geometry $\Gamma$ over $I \setminus J$ is said to be $J$-{\it locally truncated of type $D$} (or {\it with diagram $D$}) over
$I$ if and only if for every nonempty flag $F$ of $\Gamma$, the residue $\res{\Gamma}{F}$ is isomorphic to the truncation of type $I \setminus ( J \cup t(F))$ of a geometry belonging to the diagram $D_{I \setminus t(F)}$, the restriction of $D$
to the typeset $I \setminus t(F)$. If $\Gamma$ is the truncation of type $I \setminus J$ of a geometry
$\Delta$ of type $D$ over $I$ then $\Gamma$ it is a geometry of $J$-locally
truncated type $D$. The converse is in general not true; see Brouwer and Cohen
\cite{loc} and Ronan \cite{ron}.\\

Let $M=(m_{ij})$ be a Coxeter matrix with rows and columns indexed by $I$. The {\it
diagram} of $M$, denoted $D(M)$, assigns to each $2$-subset $\lbrace i, j \rbrace$ of $I$, the class $D(i,j)$ of generalized $m_{i j }$-gons. If $\Gamma$ is a residually connected geometry with diagram $D(M)$, that is every residue of type $\lbrace i, j \rbrace$ of $\Gamma$ is a generalized $m_{i j}$-gon, then $\Gamma$ is called a {\it geometry of type} $M$.

\subsection{Chamber systems} A {\it chamber system} ${\mathcal C} = (C, E, \lambda, I)$ over $I$ is a simple graph $(C, E)$ together with an edge-labeling $\lambda: E \rightarrow 2^{I} \setminus \lbrace \emptyset \rbrace$ by nonempty subsets of $I$ such that, if $a, b, c \in C$ are three pairwise adjacent vertices, then $\lambda (a,b) \cap \lambda (b,
c) \subseteq \lambda (a,c)$. The elements of $C$ are called
{\it chambers}. Two distinct chambers $a$ and $b$ are $i$-{\it adjacent} iff $(a,b) \in E$ is an edge and $i \in \lambda (a,b)$, for any $i \in I$. The rank of ${\mathcal C}$ is the cardinality of the index set $I$.\\

For $J$ a subset of $I$, the {\it residue of ${\mathcal C}$ of type $J$} or the $J$-{\it residue}, is a connected component of the graph $(C, E_J, \lambda _J, J)$, with
$\lambda_J$ the restriction of $\lambda$ to $\lambda^{-1}(2^J) \subseteq E$, where $2^J$ is the codomain of $\lambda _J$, and $E_J =\lbrace e \in E \mid \lambda _J (e) \not=
\emptyset \rbrace $. A $J$-residue $R$ is a chamber system over the typeset $J$. The set $I \setminus J$ is called the {\it cotype} of $R$. The {\it rank of the residue} $R$ is $| J|$ and its {\it corank} is $| I \setminus J |$.\\

A {\it chamber system ${\mathcal C}$ over $I$ belongs to the
diagram $D(M)$} if and only if every residue of ${\mathcal C}$ of
type $\lbrace i, j \rbrace \subseteq I$ is the chamber system of a generalized
$m_{ij}$-gon.\\

A chamber system over $I$ is {\it residually connected} if and only if for every subset $J \subseteq I$ and for every family $\lbrace R_j: j\in J \rbrace$ of residues of cotype $j$, with the property that any two have nonempty intersection, it follows that $\cap _{j \in J}R_j$ is a nonempty residue of type $I \setminus J$.\\

The chamber systems over $I$ together with the appropriate morphisms form a category denoted ${\mathcal Chamb}_I$.\\

Let $\Gamma$ be a geometry over $I$. Denote by ${\bf C}(\Gamma)$ the set of its chamber flags, that is, the flags of type $I$. Two chamber flags $c$ and $d$ are said to be $i$-adjacent whenever they have the same element of type $j$ for all $j \not = i$, with $i, j \in I$. Then ${\bf C}(\Gamma)$, with the above adjacency relation, is a chamber system of type $I$.\\

Starting with a chamber system ${\mathcal C}$ over $I$, we define a geometry ${\bf G}({\mathcal C}) = (C_i, i\in I, *, t)$. The objects of this geometry are the elements $C_i, i \in I$, the collection of all corank one residues of type $I \setminus \lbrace i \rbrace$ of ${\mathcal C}$. Two objects are incident if they have nonempty intersection.\\

The above construction gives rise to a pair of functors ${\mathbf G}: {\mathcal Chamb}_I \rightarrow
{\mathcal Geom}_I$ and
${\mathbf C}: {\mathcal Geom}_I \rightarrow {\mathcal Chamb}_I$ which have
the properties:
\begin{align*}
{\mathbf G}( {\bf C}(\Gamma)) \; &=
\; \Gamma,\; \; \; \text{if} \;  \Gamma \; \text{is a residually connected geometry and} \; I\; \text{is finite};\\
{\mathbf C}({\mathbf G}(C))\; &= \; C, \; \; \; \text{if} \; C
\; \text{is a residually connected chamber system}.
\end{align*}
For more details see Shult \cite{asp}, \cite{ln}.\\

For a connected chamber system ${\mathcal C}$ over $I$, a {\it $2$-cover} of ${\mathcal C}$ is a
connected chamber system $\widetilde {\mathcal C}$ together with a chamber system
morphism $h : \widetilde {\mathcal C} \rightarrow {\mathcal C}$ which is surjective on chambers and is an isomorphism when restricted to any residues of rank at most $2$ of $\widetilde {\mathcal C}$. A $2$-cover $h: \widetilde {\mathcal C} \rightarrow {\mathcal C}$ is
said to be {\it universal} if for any other $2$-cover $\varphi : {\mathcal C}'
\rightarrow {\mathcal C}$ there is a $2$-cover $\psi : \widetilde {\mathcal C} \rightarrow
{\mathcal C}'$ such that $h \circ \psi = \varphi$. It can be proved that
chamber systems always have universal coverings.\\

For the notion of building see the famous paper of J. Tits: {\it A local approach to buildings}, from which we reproduce the following (Corollary $3$, \cite{tits}):
\newtheorem*{LA}{Tits' Local Approach Theorem}
\begin{LA} Suppose ${\mathcal C}$ is a chamber system of type $D(M)$, $M$ a Coxeter
matrix, and suppose that for every rank $3$ residue, the universal $2$-cover is a building. Then the universal $2$-cover of ${\mathcal C}$ is a building ${\mathcal B}$ of type $D(M)$.
\end{LA}
In particular, the chamber system ${\mathcal C}$ is obtained from ${\mathcal B}$ by factoring out a group
of automorphisms in which no non-trivial element fixes any rank $2$
residue of ${\mathcal B}$.

\subsection{Sheaves}
Let $I$ be an index set, $J \subset I$ and set $K = I \setminus J$. Let $\Gamma$ be a geometry over $K$ which is locally truncated of type $D$ over $I$ and let ${\mathcal F}$ be a family of nonempty flags of $\Gamma$. A {\it sheaf over the geometry} $\Gamma$ is a class of geometries $\lbrace \Sigma(F) \; \text{for} \; F \in {\mathcal F} \rbrace$ together with isomorphisms
\begin{center}
$\varphi_F: \res{\Gamma}{F} \rightarrow \; _J \Sigma (F)$
\end{center}
of geometries over $I \setminus t(F)$. Given a pair of incident flags $F_1 \subseteq F_2$ in ${\mathcal F}$ the connecting homomorphisms of the sheaf are the maps $\varphi _{F_1, F_2} : \Sigma (F_2) \rightarrow \Sigma (F_1)$ with the property that
\begin{center} $\varphi_{F_1, F_2}(\Sigma(F_2)) \simeq \res{\Sigma(F_1)}{F_2 \setminus F_1}$.
\end{center}
Furthermore, they are subject to the following conditions:
\begin{center}
$\varphi _{F_1, F_2} \circ \varphi _{F_2, F_3} = \varphi _{F_1, F_3}$,
\end{center}
for $F_1, F_2, F_3 \in {\mathcal F}$ with $F_1 \subseteq F_2 \subseteq F_3$. To simplify the notation, we will omit the connecting isomorphisms $\varphi _F$ and write $\res{\Gamma}{F} = \; _J \Sigma(F)$ instead.\\

A sheaf $\Sigma$ is {\it residually connected} if and only if for each
object $x$ of the geometry $\Gamma$ the sheaf geometry $\Sigma (x)$ is
residually connected.\\

Due to the functorial relation between the category of geometries and the
category of chamber systems, whenever a sheaf $\Sigma$ exists, there
is also chamber system associated to it, Brouwer and Cohen (Lemma $1$ in \cite{loc}).

\subsection{Some properties of the Grassmann spaces ${\bf A_{n-1,j}}({\mathbb K})$}
In the sequel we review some of the properties of the Grassmann
spaces $A_{n-1,j}({\mathbb K})$. For details see Cohen \cite{coh}. Let {\pl} be the Grassmann space $A_{n-1,j}({\mathbb K})$. Let $V$ be a $n$-dimensional
vector space over some division ring ${\mathbb K}$. The points of $\Gamma$ are the
$j$-dimensional subspaces of $V$, the lines are the $(j-1,j+1)$-flags of $V$.
The point-line geometry $\Gamma$ is a strong parapolar space whose
symplecta are nondegenerate polar spaces of type $D_3$. Let ${\mathcal S}$ denote the family of symplecta. There are two
classes of maximal singular subspaces: ${\mathcal A}$, a
family of subspaces of singular rank $j$, and ${\mathcal B}$ whose elements have singular rank $n-j$. Set ${\mathcal M} = {\mathcal A} \cup {\mathcal B}$.

\begin{prop} Let {\pl} be the Grassmann space $A_{n-1, j}({\mathbb K})$. Then the following are true:
\begin{itemize}
\item[1.] If $S \in {\mathcal S}$ and $x \in {\mathcal P} \setminus S$ then $x^{\perp} \cap S$ is empty, a point or a plane.
\item[2.] If $S \in {\mathcal S}$ and $M \in {\mathcal M}$ then $S \cap M$ is
empty, a point or a plane.
\item[3.] If $M$ and $M'$ are two distinct maximal singular subspaces belonging to the same class then
$M \cap M'$ is either empty or a point.
\item[4.] If $M$ and $M'$ are maximal singular subspaces from different classes then $M \cap
M'$ is either empty or a line.
\end{itemize}
\end{prop}
In \cite{rem1}, Shult gave a geometric characterization, in terms of
points and maximal singular subspaces,  of
certain parapolar spaces, including some Grassmannians. We use his result here, in
order to formulate the following property of the Grassmann spaces:\\
(*) For any maximal singular subspace $M \in {\mathcal M}$
and point $x \in {\mathcal P} \setminus M$ the set $x^{\perp} \cap M$ is either empty or
a line.\\

We conclude this brief review of the Grassmannians with two technical lemmas:
\begin{lem} Let $({\mathcal P}, {\mathcal L}) \simeq A_{n-1, 4}({\mathbb K})$ with $n \geq 6$ and ${\mathbb K}$ a division ring. Take  $S \in {\mathcal S}$. Let $x \in {\mathcal P} \setminus S$ with the property that $d(x,y)=2$ for any point $y \in S$. Set ${\mathbb X}(S) = \lbrace X \in {\mathbb X} \mid X \cap S \; \text{is a plane}\; \rbrace$ where ${\mathbb X} \in \lbrace {\mathcal A}, {\mathcal B} \rbrace$. Then either $x^{\perp} \cap X = \emptyset$ for all $X \in {\mathbb X}(S)$ or $x^{\perp} \cap X$ is a line for all $X \in {\mathbb X}(S)$.
\end{lem}

\begin{proof} Consider a vector space $V$ of dimension $n \geq 6$ (all the dimensions are affine dimensions) over some division ring ${\mathbb K}$. The points ${\mathcal P}$ are the $4$-subspaces of $V$. Take a $2 - 6$ flag $F_2 \subset F_6$ in $V$ and let $S$ be the symplecton determined by this flag, that is, all the $4$-subspaces containing $F_2$ and contained in $F_6$. Recall $x$ is a $4$-subspace such that, for any $4$-subspace $y$ with $F_2 \subset y \subset F_6$, the set $x \cap y$ is a $2$-subspace. Three cases are possible: $a)\; \; \text{dim}(x \cap F_2) = 2;\;\; b)\; \; \text{dim}(x \cap F_2) = 0; \; \; c) \; \; \text{dim}(x \cap F_2) = 1$.

\medskip
We start by eliminating case $c)$. So assume $x$ is such that $x \cap F_2$ is a $1$-dimensional subspace. Note that $\text{dim}(x \cap F_6)$ can be at most $4$. If $\text{dim}(x \cap F_6) =1$ then $\text{dim}(x \cap y) \leq 1$, for any $4$-subspace $y \in S$, which implies $d(x, y) \geq 3$ and this contradicts the hypothesis on $x$. If $\text{dim}(x \cap F_6) = 2$, take a $2$-subspace $U \subset F_6$ to be such that $U \cap x = 0$ and $U \cap F_2 = 0$ and let $y = \langle F_2, U \rangle $. Then $\text{dim}(x \cap y) \leq1$ and $d(x, y) \geq 3$, a contradiction. Let now $\text{dim}(x \cap F_6) \in \lbrace 3, 4 \rbrace$. Let $y \in S$ be such that $y \subseteq \langle F_2, x \cap F_6 \rangle $. Then $\text{dim}(x \cap y) = 3$ which is a contradiction. Therefore $\text{dim}(x \cap F_2) \not= 1$.

\medskip
When $\text{dim}(x \cap F_2) = 2$, since $d(x, y) =2$ for any $y \in S$, then $x \cap F_6 = x \cap F_2 = F_2$.

\medskip
Let now $x \cap F_2 = 0$. In this case $x \subset F_6$. Otherwise $\text{dim}(x \cap F_6) \leq 3$ and then we can find $y \in S$ to be such that $\text{dim}(x \cap y) \leq 1$, contradicting the hypothesis on $x$. \\
We shall discuss each class of maximal singular subspaces separately.\\
First let ${\mathbb X} = {\mathcal A}$. Let $A \in {\mathcal A}(S)$, that is all $4$-subspaces contained in a fixed $5$-subspace $F_A$. Since $A \cap S$ is a plane, it follows $F_2 \subset F_A \subset F_6$. If $x$ is of type $a)$, i.e. $x \cap F_6 = F_2$ then $x \cap F_A = F_2$. Consequently $x^{\perp} \cap A = \emptyset$. If $x$ is of type $b)$, that is $x \cap F_2 = 0$ and $x \subset F_6$ then $\text{dim}(x \cap F_A)=3$. Note that $x \not \in A$ because this would imply $x^{\perp} \cap S \not = \emptyset$, contradicting the hypothesis on $x$. In this case $x^{\perp} \cap A$ contains the $4$-subspaces $y$ with the property $x \cap F_A \subset y \subset F_A$. Therefore, $x ^{\perp} \cap A$ is a line; this also follows from (*). Since $A$ is arbitrary in ${\mathcal A}(S)$, the conclusion follows in this case.\\
Let now ${\mathbb X} = {\mathcal B}$ and let $B \in {\mathcal B}$. Then $B$ is determined by all $4$-subspaces of $V$ containing a fixed $3$-subspace $F_B$. Assume $B \cap S \not = \emptyset$ and let $y \in B \cap S$. Then since $y \in S, \; F_2 \subset y \subset F_6$ and since $y \in B, \; F_B \subset y$. Therefore $\text{dim}(F_2 \cap F_B)$ is either $1$ or $2$. If $\text{dim}(F_2 \cap F_B) = 1$ then $B \cap S = \langle F_2, F_B \rangle $ is a unique point (see also Proposition $2.1$, part $2$). If $\text{dim}(F_2 \cap F_B) = 2$ then $F_2 \subset F_B \subset F_6$ and therefore in this case $B \cap S$ is a plane. Thus $B \in {\mathcal B}(S)$ is equivalent with the condition that $F_2 \subset F_B \subset F_6$.\\
Assume that $x$ is as in $a)$. Let $U \subset x$ be a $1$-subspace with $U \cap F_2 = 0$. Take $y = \langle F_B, U \rangle $. Then $y \in x^{\perp} \cap B$ and therefore (*) implies $x^{\perp} \cap B$ is a line.\\
Let now $x$ be as in $b)$. Let $y \in B$ be an arbitrary point. Then $x \cap y \subseteq y$ and $F_2 \subset y$ are two disjoint $2$-subspaces since $x \cap F_2 = 0$. Therefore $\text{dim}(x \cap y) \leq 2$, which implies $x^{\perp} \cap B = \emptyset$. Since $B$ is arbitrarily in ${\mathcal B}(S)$, the conclusion follows now.
\end{proof}

Throughout the next Lemma, the notations introduced in Lemma $2.2$ are maintained.

\begin{lem}Let $\Gamma = ({\mathcal P}, {\mathcal L}) \simeq A_{n-1, 4}({\mathbb K})$ with $n > 6$. Assume that $S, R \in {\mathcal S}$ are two symplecta such that $S \cap R$ is a plane $\pi$. Let $x \in {\mathcal P}$ be a point which is at distance two from every point of $S$. Assume that $x^{\perp} \cap X$ is a line for some $X \in {\mathbb X}(S)$. Then $x^{\perp} \cap R$ is either empty or a plane.\end{lem}

\begin{proof} Let $S$ and $R$ be two symplecta as in the hypothesis. Assume that $S$ is determined by the $2 - 6$ flag $F_2 \subset F_6$ and $R$ is determined by the flag $F'_2 \subset F'_6$. Then the points $y \in S \cap R$ are the $4$-subspaces with the property that $\langle F_2, F'_2\rangle  \subseteq y \subseteq F_6 \cap F'_6$. Now $\langle F_2, F'_2 \rangle $ is at least $2$-dimensional and $F_6 \cap F'_6$ can be at most $6$-dimensional. Therefore the points $y$ will comprise a plane if either $(i) \; \text{dim}\langle F_2, F'_2\rangle =3$ and $\text{dim}(F_6 \cap F'_6) = 6$ or $(ii)\; \text{dim}\langle F_2, F'_2\rangle  = 2$ and $\text{dim}(F_6 \cap F'_6 )= 5$.

\medskip
Assume first that $(i) \; \text{dim}\langle F_2, F'_2\rangle =3$ and $F_6 = F'_6$. Therefore $\pi = R \cap S$ is the collection of $4$-subspaces containing a fixed $3$-subspace $F_B$ and contained in the $6$-subspace $F_6 = F'_6$. Let $x \in {\mathcal P}$ be a point at distance two from every point in $S$ and such that $x^{\perp} \cap X$ is a line for some $X \in {\mathbb X}(S)$. According to the result of the previous Lemma, $x^{\perp} \cap X$ is a line for every $X \in {\mathbb X}(S)$. From the proof of the previous Lemma it follows that either $a.) \; \text{dim}(x \cap F_2) = 2$ or $b.) \; \text{dim}(x \cap F_2)=0$. If $\text{dim}(x \cap F_2) = 2$ then $x \cap F_6 = F_2$ (see the proof of Lemma $2.2$ again). Since $F_2 \subset x$ it follows now that $\text{dim}(x \cap F'_2) = \text{dim}(x \cap F'_2 \cap F_2) =1$. Also, in this case $\text{dim}(x \cap F'_6)=2$. If $y \in R$ then $x \cap y \subseteq x \cap F'_6 = x \cap F_6 = F_2$ and therefore $x \cap y$ is at most $2$-dimensional and this implies $x^{\perp} \cap R = \emptyset$.\\
Assume now that $b.)\; x \cap F_2 = 0$. Using the previous Lemma, we conclude that $x \subset F_6$ in this case. Therefore $\text{dim}(x \cap F'_2) \leq 1$. If $\text{dim}(x \cap F'_2) = 0$ then $x^{\perp} \cap R = 0$. Let now $\text{dim}(x \cap F'_2) = 1$. The points $ y \in x^{\perp} \cap R$ are precisely the $4$-subspaces $F'_2 \subset y \subset \langle F'_2, x \rangle $ and this collection is a plane.

\medskip
The proof in the second case is dual to part $(i)$, where intersections and spans are interchanged. Let now $(ii) \; \text{dim}\langle F_2, F'_2\rangle  = 2$ and $\text{dim}(F_6 \cap F'_6 )= 5$. In this case $\pi = S \cap R$ can be described as the collection of the $4$-subspaces containing a $2$-subspace $F_2 = F'_2$ and contained in a fixed $5$-subspace $F_A = F_6 \cap F'_6$. Let now $x \in \mathcal{P}$ be a point at distance two from every point of $S$ and such that $x^{\perp} \cap X$ is a line for some $X \in {\mathbb X}(S)$. From the previous Lemma it follows that $x ^{\perp} \cap X$ is a line for every $X \in {\mathbb X}(S)$. Again there are two cases to consider: $a.) \; \text{dim}(x \cap F_2) = 2$ and $b.) \; \text{dim}(x \cap F_2) = 0$. Let us first assume that $a.)\; \text{dim}(x \cap F_2) = 2$. Therefore $x$ is such that $x \cap F_6=F_2$ and $x \cap F_A = F_2$. We study in more detail the relation between $x$ and $F'_6$. Clearly $\text{dim}(x \cap F'_6)$ can be $2, 3, 4$. If $\text{dim}(x \cap F'_6)=4$, since $x \cap F'_6$ and $F_A$ are both subspaces of $F'_6$ it follows $\text{dim}(x\cap F'_6 \cap F_A) \geq 3$ contradicting the properties of $x$. If $\text{dim}(x \cap F'_6)=2$ then, given any $y \in R$, $\text{dim}(x \cap y) \leq 2$ and in this case $x^{\perp} \cap R = \emptyset$. We are left with the case $\text{dim}(x \cap F'_6) = 3$. Since $F_2 = F'_2$ then $F'_2 \subset x$ and consequently all the $4$-subspaces of $F'_6$ containing the $3$-subspace $x \cap F'_6$ are points $y \in x^{\perp} \cap R$. In this case $x^{\perp} \cap R$ is a plane.\\
Let now $b.) \; \text{dim}(x \cap F_2) = 0$. In this case $x \subset F_6$ and $\text{dim}(x \cap F_A) = 3$; see the proof of Lemma $2.2$. Therefore $x \cap F'_2 = 0$ and $\text{dim}(x \cap F'_6)$ can be $3$ or $4$. In both cases $x^{\perp} \cap R = \emptyset$. This concludes the proof of the Lemma.
\end{proof}

\subsection{Some properties of ${\bf D_{5,5}}({\mathbb K})$ and ${\bf D_{6,6}}({\mathbb K})$}
The half-spin geometry $D_{n,n}({\mathbb K})$ is the point-line geometry whose points are the
maximal singular subspaces belonging to one of the classes of the oriflamme polar
space $D_n({\mathbb K})$ and the lines are the singular subspaces of dimension $n-2$ of the
same oriflamme geometry. ${\mathbb K}$ is a field. The half-spin geometries are
strong parapolar spaces whose
symplecta are nondegenerate polar spaces of type $D_4$. Let ${\mathcal S}$ denote the family of symplecta. There are two classes of
maximal singular subspaces: ${\mathcal A}$, whose elements have singular rank $n-1$, and ${\mathcal B}$, whose elements have singular rank $3$. The half-spin geometries for $n=5,6$ were characterized by Cohen and Cooperstein \cite {cc}:
\newtheorem*{HS}{Theorem}
\begin{HS}(Cohen and Cooperstein). Let {\pl} be a strong parapolar space, not
a polar space, whose symplecta have polar rank $4$. Assume that, given a
non-incident point-symplecton pair $(x, S) \in {\mathcal P} \times {\mathcal S}$ the set $x^{\perp} \cap S$ is either a point or a maximal singular subspace of $S$. Then
$\Gamma \simeq D_{6,6}({\mathbb K})$ or $\Gamma \simeq D_{5,5}({\mathbb K})$.
\end{HS}

\section{The geometry $\Gamma$}
In what follows {\pl} will be a parapolar space which is locally $A_{n-1,j}({\mathbb K})$, where ${\mathbb K}$ is a field, $n$ is an integer greater than $6$ and $j$ is $3$ or $4$. This means that at every point $p \in {\mathcal P}$, the geometry of lines and planes on $p$, denoted $\res{\Gamma}{p}$, is the geometry $A_{n-1, j}({\mathbb K})$. Let ${\mathcal S}$ denote the set of symplecta in $\Gamma$, these are nondegenerate polar spaces of type $D_4$.\\

{\it Notation:} Let $X$ be a subspace of {\pl} and assume that ${\mathcal F}$ is a family of subspaces of $\Gamma$, which does not necessarily contain $X$. Then the set of all the elements $F \in {\mathcal F}$ which are incident with $X$ will be denoted ${\mathcal F}(X)$. To simplify the notation, when $X = \lbrace p \rbrace$ is a point, we shall write ${\mathcal F}_p = {\mathcal F}(\lbrace p \rbrace)$ for the collection of those elements $F$ of ${\mathcal F}$ which contain $p$. \\

From the local assumption, it follows that, given a point $p \in {\mathcal P}$, the set  ${\mathcal M}_p$ of the maximal singular subspaces at $p$ partitions in two classes ${\mathcal A}_p$ and ${\mathcal B}_p$. If $n \not = 2j$ then the above partition is realized according to the dimension. The point $p$ was arbitrarily chosen and, given any other point $q \in {\mathcal P}$, the residual geometries have the same type, that is $\res{\Gamma}{p} \simeq \res{\Gamma}{q}$. Therefore we recognize a global partition of the maximal singular subspaces ${\mathcal M} = {\mathcal A} \cup {\mathcal B}$. If $n = 2j$ then the elements of ${\mathcal A}_p$ have the same dimension as the elements of ${\mathcal B}_p$. A global partition of the maximal singular subspaces cannot be obtained anymore (details on this special case can be found in Onofrei \cite{thesis}).\\

For the rest of this paper we shall assume that $n \not = 2j$. Therefore $\Gamma$ has the following properties, which are direct consequences of the theory of parapolar spaces and of the properties of Grassmann spaces listed in Subsection $2.5$: \\

{\bf (L.1)} There are two classes of maximal singular subspaces
${\mathcal A}$ and ${\mathcal B}$.  Set ${\mathcal M} = {\mathcal
A} \cup
{\mathcal B}$. Two distinct maximal singular
subspaces which belong to the same class meet at a line, a
point or the empty set. Two maximal singular subspaces from different
classes can meet at a plane, a point or they are disjoint.\\

{\bf (L.2)} Given $(S, M) \in {\mathcal S} \times {\mathcal M}$ then $S
\cap M$ is empty, a point, a line  or a
maximal singular subspace of $S$.

\medskip
{\bf (L.3)} For $S \in {\mathcal S}$ let ${\mathcal  M}(S)$ be the family of maximal singular subspaces with largest intersection with $S$. It is a direct consequence of $(L.1)$ that ${\mathcal M}(S) = {\mathcal A}(S) \cup {\mathcal B}(S)$. Let $M_i \in {\mathcal M}(S),\;  i=1,2$. Denote by
${\overline M}_i = M_i \cap S$. Then ${\overline M}_1 \cap {\overline M}_2$ is a
point or a plane if
${\overline M}_1, {\overline M}_2$ belong to different classes. If they belong to the same class then
${\overline M}_1 \cap {\overline M}_2$ is empty, a line or they are equal.

\medskip
{\bf (L.4)} Given $M \in {\mathcal M}$ and $x \in
{\mathcal P} \setminus  M$ then
$x^{\perp} \cap M$ is empty, a point or a plane.

\medskip
{\bf (L.5)} Given $S \in {\mathcal S}$, $x \in {\mathcal P} \setminus S$
then $x^{\perp} \cap S$ is empty, a point, a line or a maximal singular subspace of
$S$.\\

For $S \in {\mathcal S}, \; {\mathbb X} \in \lbrace {\mathcal A}, {\mathcal B}
\rbrace$ and $X \in \lbrace A, B \rbrace$ define the following sets:
\begin{center}
$M_{\mathbb X}(S)= \lbrace {\overline X} \mid {\overline X} = X \cap  S
\; \; \text{for some} \; X \in {\mathbb X}(S) \rbrace$
\end{center}
which are the two classes of maximal singular subspaces of $S$.
Then set:
\begin{center}
$N_{\mathbb X}(S)= \lbrace x \in {\mathcal P} \setminus S \mid
x^{\perp}
\cap S \in M_{\mathbb X}(S) \rbrace$
\end{center}
Also define: \begin{equation*}
\begin{split}
{\mathcal X}(S) & = \lbrace x \in {\mathcal P} \setminus S \mid
x^{\perp} \cap S =
\lbrace p \rbrace, \; p \in {\mathcal P};\; \text{for any}\; q \in p^{\perp}\cap S\\
&\quad \text{the pair} \; \lbrace x,q \rbrace \; \text{is polar}; \; \text{for some} \; A \in {\mathcal A}_p \cap {\mathcal A}(S), \; x^{\perp} \cap A \; \text{is a plane} \rbrace
\end{split}
\end{equation*}

{\it Notation:} In what follows if $x \in {\mathcal
P},\ S \in {\mathcal S}$ are such that  $x ^{\perp} \cap S \in M_{\mathbb X}(S)$,
then we denote $x^{\perp} \cap S = {\overline X}_x$ and the maximal singular
subspace containing it $X_x$. Next, if $p, q \in {\mathcal P}$ form a polar pair in $\Gamma$, the unique symplecton containing them will be denoted $\ll x, y \gg$.

\begin{lem} Let $x \in N_{\mathbb X}(S)$ and $y \in x^{\perp}$ be such that
$y^{\perp} \cap S \setminus {\overline X}_x \not= \emptyset$. Then $y \in S \cup N_{\mathbb X}(S)$.
\end{lem}
\begin{proof} If $y \in S$ we are done. Assume that $y \not\in S$ and let $p \in y^{\perp} \cap S \setminus {\overline X}_x$, which exists by hypothesis. Set $R= \ll x,p \gg$ and note that $y \in R$.
Therefore $y^{\perp} \cap p^{\perp} \cap \overline X_x =  L$, a line. Then $y^{\perp} \cap S$ contains the plane $\langle p, L \rangle $. According to $(L.5)$, $y^{\perp} \cap S = {\overline M}$ is a maximal singular subspace of $S$. In order to prove that $y \in N_{\mathbb X}(S)$ it suffices
to show that ${\overline M} \cap {\overline X}_x = L$. Clearly ${\overline M} \cap {\overline X}_x \supseteq L$. Assume by contradiction that ${\overline M} \cap {\overline X}_x$ contains a plane $\pi$. Then $\langle xy, \pi \rangle $ is a singular subspace of rank at least $4$ inside
the symplecton $R$. This is a contradiction since $R$ is a rank $4$ polar space.
\end{proof}

Throughout the rest of the paper we shall occasionally use {\it the gamma space property of $\Gamma$}. Recall that this means that, given $(x, L) \in {\mathcal P} \times {\mathcal L}$ with $x \not \in L$, then $x^{\perp} \cap L$ can be empty, a point or the entire line $L$.

\begin{lem} For any symplecton $S \in {\mathcal S}$ the set $S \cup N_{\mathbb X}(S)$ is a
subspace of $\Gamma$.
\end{lem}
\begin{proof}
Let $x$ and $y$ be two collinear points in $S \cup N_{\mathbb X}(S)$.
If at least one of the points is in $S$ the conclusion follows at
once. So we may assume that $x, y \not\in S$. Then
${\overline X}_x = x^{\perp} \cap S$ and ${\overline X}_y = y^{\perp} \cap S$ are two
maximal singular subspaces of $S$ belonging to the same class. According to $(L.3)$
there are three cases to consider:

\medskip
  {\bf (i).} If ${\overline X}_x = {\overline X}_y$ then, by $(L.1)$, $X_x = X_y$ and
$xy \subset N_{\mathbb X}(S)$.

\medskip
 {\bf (ii).} Assume now that ${\overline X}_x \cap {\overline X}_y = L$, a line in $S$.
Let $z \in xy \setminus \lbrace x, y \rbrace$. By the gamma space property $L \subset z^{\perp}$. Let $w \in
{\overline X}_y \setminus z^{\perp}$. Now ${\overline N}= \langle w, w^{\perp} \cap {\overline X}_x \rangle $ is a maximal singular subspace of $S$, not in the
same class with ${\overline X}_x$ and ${\overline X}_y$. This is
true since ${\overline N} \cap {\overline X}_x$ and ${\overline N} \cap {\overline
X}_y$ are both planes in $S$. Set $R =  \ll w,x \gg$ and notice that $\overline N \subseteq R \cap S$. In $R$,
$z^{\perp} \cap {\overline N}$ is a plane and therefore, by $(L.5)$, $z ^{\perp}
\cap S = {\overline X}_z$ is a maximal singular subspace of $S$. It remains to show that
${\overline X}_z \in \mss{\mathcal X}{S}$. It suffices to prove: ${\overline X}_x \cap {\overline X}_z =L$. Now ${\overline X}_z$ meets ${\overline N}$ at a plane, namely the plane $z^{\perp} \cap {\overline N}$. So, the family of ${\overline X}_z$ is not the same as that of ${\overline N}$, which in its turn is not the same as that of ${\overline X}_x$ and ${\overline X}_y$. So ${\overline X}_z$ and ${\overline X}_x$, being different, meet at a line, which must be $L$.

\medskip
{\bf (iii).} Assume that ${\overline X}_x \cap {\overline X}_y = \emptyset$. Let $u
\in {\overline X}_x$ and set $R =  \ll u, y \gg$. Then the plane $u^{\perp} \cap {\overline X}_y$ lies in $R$ and therefore $x^{\perp} \cap u^{\perp} \cap {\overline X}_y$ is a line. But this is a contradiction with the assumption that ${\overline X}_x \cap {\overline X}_y = \emptyset$. Thus
either ${\overline X}_x \cap {\overline X}_y \not= \emptyset$ or $x$ is not collinear with $y$.
\end{proof}

\begin{lem} Let $S \in {\mathcal S}$ and $x, y \in N_{\mathcal A}(S)$ be
two points at distance two. Assume that ${\overline A}_x \cap {\overline A}_y = L$, a line. Let $z \in x^{\perp} \cap y^{\perp} \setminus S$ be such that $L \subset z^{\perp}$. Then $z \in A_x \cup A_y$.
\end{lem}
\begin{proof} Set $R = \ll x, y \gg$. By
hypothesis $L \subset z^{\perp}$ so $\langle xz, L\rangle $ and $\langle zy, L\rangle $ are two maximal singular subspaces in $R$. By $(L.3)$, they belong to different
classes in $R$. Therefore either $\langle xz, L\rangle  \subseteq A_x$ or $\langle zy,L\rangle  \subseteq A_y$. Consequently $z \in A_x \cup A_y$.
\end{proof}

\bigskip
\section{The class ${\mathbb D}$ of subspaces}
Throughout this section {\pl} is a parapolar space, locally $A_{n-1,3}({\mathbb K})$
for some field ${\mathbb K}$ and $n$ an integer greater than $6$. The properties $(L.1)$ to $(L.5)$ listed at
the beginning of the previous section remain valid and the notations introduced there are maintained. In this case the maximal singular subspaces from class ${\mathcal A}$ have singular rank $4$. Set ${\mathsf D}(S) = S \cup \nas{S}$ and let ${\mathbb D} = \lbrace {\mathsf D}(S) \mid S \in {\mathcal S} \rbrace$. We prove the following:
\begin{ta} Let {\pl} be a parapolar space which is locally $A_{n-1,3}({\mathbb K})$
for some integer $n > 6$ and field ${\mathbb K}$. There is a
collection ${\mathbb D}$ of $2$-convex subspaces, each isomorphic to
$D_{5,5}({\mathbb K})$ such that, if $S \in {\mathcal S}$ is a symplecton, then
there exists a unique member ${\mathsf D}(S) \in {\mathbb D}$ containing $S$.
\end{ta}

Recall that by Lemma $3.2, \; {\mathsf D}(S)$ is a subspace of $\Gamma$.

\begin{lem} Let $R, S \in {\mathcal S}$. If $R \cap S \in \mss{\mathcal B}{S}$ then
${\mathsf D}(S) = {\mathsf D}(R)$.
\end{lem}

\begin{proof} Let $S, R \in {\mathcal S}$ be such that $S \cap R = {\overline B}$ is a maximal singular subspace of ${\mathcal B}$-type of the two symplecta. It will suffice to prove that ${\mathsf D}(S) \subseteq {\mathsf D}(R)$.\\
Let $x$ be a point in ${\mathsf D}(S)$. If $x \in {\overline B} \subset R$ we are done. Assume now that $x^{\perp} \cap {\overline B}$ is a plane. This is true if either $x \in S \setminus {\overline B}$ or $x \in \nas{S}$ and $x^{\perp} \cap S \cap {\overline B}$ is a plane. Then, by $(L.5)$, $x^{\perp} \cap R$ has to be a maximal singular subspace of $R$. Set $x^{\perp} \cap R = {\overline A'}_x$. Since ${\overline A'}_x \cap {\overline B}$ is a plane, it follows, by $(L.3)$, that ${\overline A'}_x$ and ${\overline B}$ belong to different families of maximal singular subspaces of $R$. Therefore ${\overline A'}_x \in \mss{\mathcal A}{R}$ and $x \in \nas{R}$.

\medskip
Let now $x \in \nas{S}$ and set ${\overline A}_x= x^{\perp} \cap S$. Assume that ${\overline A}_x \cap {\overline B} = \lbrace p \rbrace$, a single point. Let $q \in {\overline B} \setminus \lbrace p \rbrace $ and note that $\lbrace x, q \rbrace$ is a polar pair, since $x^{\perp} \cap q^{\perp}$ contains the plane $q^{\perp} \cap {\overline A}_x$. So we can find a point $w \in p^{\perp} \cap q^{\perp} \cap x^{\perp} \setminus S$. By Lemma $3.1$, applied to the pair $\lbrace x, w \rbrace$ and symplecton $S,\; w \in \nas{S}$. According to the previous paragraph $w \in {\mathsf D}(R)$. Let now $y \in {\overline A}_x \setminus w^{\perp}$ so $y \in \nas{R}$. If $w \in R$, apply Lemma $3.1$ to the pair $\lbrace x, y \rbrace$ and symplecton $R$ and get $x \in \nas{R}$. So let us assume now that $w \not \in R$. Let $A_w$ be the maximal singular subspace $\langle w, w^{\perp} \cap R \rangle $. Recall that $A_w$ has singular rank $4$. Now $x^{\perp} \cap A_w$ contains the line $pw$ and, by $(L.4)$, $x^{\perp} \cap A_w$ is a plane. Then $x^{\perp} \cap A_w$, a plane, and $A_w \cap R$, which has singular rank $3$, must have at least a line in common. It follows that $x^{\perp} \cap R$ contains a line, say $px_1$, with $x_1 \in R \setminus {\overline B}$. Now $x_1 \not \in y^{\perp}$ since $x_1^{\perp} \cap S = A_w \cap S$ and $y \not \in w^{\perp}$. Since $y \in \nas{R}$, we may apply Lemma $3.1$ to the pair of points $\lbrace x, y \rbrace$ and symplecton $R$ and conclude that $x \in \nas{R}$. This finishes the proof that ${\mathsf D}(S) \subseteq {\mathsf D}(R)$.
\end{proof}

\begin{prop} For $S \in {\mathcal S},\ {\mathsf D}(S)$ is a $2$-convex subspace of $\Gamma$.
\end{prop}

\begin{proof} Let $x, y \in {\mathsf D}(S)$ be two points at distance two. We have to prove that $x^{\perp} \cap y^{\perp} \subset {\mathsf D}(S)$. If $x, y \in S$ then we are done since $S$ is already a $2$-convex subspace of $\Gamma$.

\medskip
Let $x \in \nas{S}$, with ${\overline A}_x = x^{\perp} \cap S$, and $y \in S$. Then $\lbrace x, y \rbrace$ is a polar pair. Set $R= \ll x,y \gg$. Now $R \cap S = \langle y, y^{\perp} \cap {\overline A}_x \rangle  \in \mss{\mathcal B}{S}$. Then, by Lemma $4.1, \; {\mathsf D}(S) = {\mathsf D}(R)$, which implies that $x^{\perp} \cap y^{\perp} \subset {\mathsf D}(S)$.

\medskip
Let now $x, y \in \nas{S}$ and set ${\overline A}_x = x^{\perp} \cap S$ and ${\overline A}_y = y^{\perp} \cap S$. Let $p \in {\overline A}_x \setminus {\overline A}_y$. Set $R = \ll p, y \gg$. Then $S \cap R = \langle p, p^{\perp} \cap {\overline A}_y \rangle  \in M_{\mathcal B}(S)$. Therefore, by Lemma $4.1, \; {\mathsf D}(S) = {\mathsf D}(R)$. Now $x \in {\mathsf D}(R), \; y \in R$ and according to the previous steps of this proof, $x^{\perp} \cap y^{\perp} \subset {\mathsf D}(R) = {\mathsf D}(S)$.
\end{proof}

\begin{prop} Given $S \in {\mathcal S}, \; {\mathsf D}(S)$ is a strong parapolar subspace of $\Gamma$.
\end{prop}

\begin{proof} We have to prove that, given two points $x, y \in {\mathsf D}(S)$ at distance two, $x^{\perp} \cap y^{\perp}$ contains at least two points.

\medskip
  {\bf (i).} If $x, y \in S$ we are done since $\ll x, y \gg =S$.

\medskip
 {\bf (ii).} Let $x \in \nas{S}$ and $y \in S$. Then $x^{\perp} \cap y^{\perp} \cap S$ is a plane and therefore the pair $\lbrace x, y \rbrace$ is polar.

\medskip
{\bf (iii).} If $x, y \in \nas{S}$ then ${\overline A}_x = x^{\perp} \cap S$ and ${\overline A}_y = y^{\perp} \cap S$ can meet at a line or the empty set. If ${\overline A}_x \cap {\overline A}_y = L$, a line, then $x^{\perp} \cap y^{\perp} \supset L$ and $\lbrace x, y \rbrace$ is a polar pair. So assume that ${\overline A}_x \cap {\overline A}_y = \emptyset$. Let $p \in {\overline A}_x$ and set $R = \ll p, y \gg$. Now $R \cap S = \langle p, p^{\perp} \cap {\overline A}_y \rangle  \in \mss{\mathcal B}{S}$ and by Lemma $4.1, \; {\mathsf D}(R) = {\mathsf D}(S)$. Then $x \in \nas{S} \subset {\mathsf D}(R)$. Since $y \in R$, by Steps $(i)$ and $(ii)$ of this Proposition, it follows that $\lbrace x, y \rbrace$ is a polar pair.
\end{proof}

\begin{prop} Let $S \in {\mathcal S}$ and let $R$ be some symplecton in ${\mathsf
D}(S)$. Then ${\mathsf D}(S) = {\mathsf D}(R)$. That is, any
symplecton $R \in {\mathcal S}$ is contained in a unique ${\mathsf D}(S)$
for some $S \in {\mathcal S}$.
\end{prop}
\begin{proof}  Let $S, R \in {\mathcal S}$ be such that $R \subset
{\mathsf D}(S)$ and $R \cap S
\not= \emptyset$. Let $ p \in R \cap S$ and $ x \in R \setminus p^{\perp}$. Since $R \subset {\mathsf D}(S)$, $ x ^{\perp} \cap S =
{\overline A}_x$ a maximal singular subspace of $S$. Then the plane $x^{\perp} \cap p^{\perp} \cap S$ lies in $\ll p, x \gg = R$ and also in $S$. Therefore $R \cap S = \langle p, p^{\perp} \cap {\overline A}_x \rangle  \in \mss{\mathcal B}{R} \cap \mss{\mathcal B}{S}$. Now, according to Lemma $4.1$, ${\mathsf D}(R) = {\mathsf D}(S)$.

\medskip
Let us now assume that $S \cap R = \emptyset$ and $R \subset {\mathsf
D}(S)$. Let $x \in R$, then $x ^{\perp} \cap S = {\overline A}_x \in \mss {\mathcal A}{S}$. Let $y
\in S \setminus {\overline A}_x$. Then set $T =  \ll x, y \gg$. Now $T \cap S = \langle y, y^{\perp} \cap {\overline A}_x\rangle
\in \mss{\mathcal B}{S}$ and, by Lemma $4.1,\ {\mathsf D}(S) = {\mathsf D}(T)$.
Also $R \cap T \not= \emptyset$ since $x \in R \cap T$ and, by the result of the previous paragraph,
it follows that ${\mathsf D}(R) = {\mathsf D}(T) = {\mathsf D}(S)$.
\end{proof}
\bigskip
\begin{proof}[{\bf Proof of the Theorem 1}] Let {\pl} be a parapolar space
which is locally $A_{n-1,3}({\mathbb K})$ for $n > 6$ and ${\mathbb K}$ a field. Given $S \in {\mathcal S}$, a symplecton, there is a $2$-convex subspace ${\mathsf D}(S)$ containing it, by Lemma $3.2$ and Proposition $4.2$. Also, by Proposition $4.3, \; {\mathsf D}(S)$ is a strong parapolar subspace of $\Gamma$. Moreover, according to Proposition $4.4$, for any symplecton $R \in {\mathcal S}$ there is a unique subspace ${\mathsf D}(S) ={\mathsf D}(R)$ containing $R$, for some $S \in {\mathcal S}$. If $x \in {\mathsf D}(S)$ is a point and $R \subset {\mathsf D}(S)$ some symplecton, since ${\mathsf D}(S) = {\mathsf D}(R)$, then either $x \in R$ or $x^{\perp} \cap R \in \mss{\mathcal A}{R}$.
Using the characterization of the half-spin geometries given by
Cohen and Cooperstein \cite{cc}, see the Theorem from Section $2.6$, we identify the subspace ${\mathsf D}(S)$ with $D_{5,5}({\mathbb K})$.
\end{proof}

\section{The class ${\mathcal D}$ of subspaces}
Throughout this section {\pl} is a parapolar space which is locally $A_{n-1,4}({\mathbb K})$, with $n =7$ or $n$ an integer greater than $8$. ${\mathbb K}$ is a field. Therefore $\Gamma$ has the properties $(L.1) - (L.5)$ listed at the beginning of Section $3$. The notations from Section $3$ are maintained. In this case the maximal singular subspaces which belong to the class ${\mathcal A}$ have singular rank $5$.
Set $D(S) = S \cup {\mathcal X}(S) \cup \nas{S}$ and let ${\mathcal D} = \lbrace D(S) \mid S \in {\mathcal S} \rbrace$.\\

Furthermore we assume $\Gamma$ satisfies the following:\\
{\bf The Weak Hexagon Axiom} {\it (WHA)}: {\it Let $H = (x_1, \ldots x_6)$ be a
$6$-circuit, isometrically embedded in $\Gamma$, this means that $x_i \in x^{\perp}_{i+1}$, indices taken mod $6$, and all the other pairs are not collinear. Also assume that at least one of the pairs of points at distance two, say $\lbrace x_1, x_3 \rbrace$, is polar. Then there exists a point $w \in x_1^{\perp} \cap x_3^{\perp} \cap x_5^{\perp}$}.\\

The purpose of this section is to prove:
\begin{tb} Let {\pl} be a parapolar space which is locally $A_{n-1,4}({\mathbb K})$
for some $n=7$ or $n \geq 9$ and field ${\mathbb K}$. Assume $\Gamma$ satisfies (WHA). Then there is a
collection ${\mathcal D}$ of $2$-convex subspaces, each isomorphic to the
half-spin geometry $D_{6,6}({\mathbb K})$ such that, if $S \in {\mathcal S}$ is a
symplecton, there exists a unique member $D(S) \in {\mathcal D}$
containing $S$.
\end{tb}
Let $S \in {\mathcal S}$ be a symplecton in $\Gamma$. We first prove that $D(S)$ is a $2$-convex strong parapolar subspace of $\Gamma$. In order to prove that $D(S) \simeq D_{6,6}({\mathbb K})$ we use Cohen-Cooperstein characterization theorem mentioned in Section $2.6$. We start with a few technical lemmas.

\begin{lem} Let $x \in {\mathcal X}(S)$ and denote $\lbrace p \rbrace = x^{\perp} \cap S$. Then, for any $A \in {\mathcal A}_p \cap {\mathcal A}(S),\; x^{\perp} \cap A$ is a plane.
\end{lem}

\begin{proof} Let $F$ be a subspace of $\Gamma$ on $p$ (since all the subspaces considered in the sequel contain the point $p$ we shall omit the subscript $p$); we shall denote by ${\widetilde F}$ the corresponding subspace of $\res{\Gamma}{p}$ whose ``points" are the lines of $F$ containing $p$ and whose ``lines" are the planes of $F$ on $p$. If $F$ belongs to a family ${\mathcal F}$ of subspaces in $\Gamma$ on $p$, ${\widetilde {\mathcal F}}$ will denote the corresponding class of subspaces in $\res{\Gamma}{p}$.\\
Since $x \in {\mathcal X}(S)$, $d_{\res{\Gamma}{p}}({\tilde x}, {\tilde q})=2$ for any ``point" ${\tilde q} \in {\widetilde S}$. Furthermore, according to the definition of ${\mathcal X}(S)$, there is a ``maximal singular subspace" ${\tilde A}_0\in {\tilde {\mathcal A}}$, of singular rank $4$, such that ${\tilde A}_0 \cap {\widetilde S}$ is a ``plane" and ${\tilde x}^{\tilde \perp} \cap {\tilde A}_0$ is a ``line". Apply Lemma $2.2$ to the ``point"-``symplecton" pair $({\tilde x}, {\tilde S})$ from $\res{\Gamma}{p}$ and conclude that ${\tilde x}^{\tilde \perp} \cap {\tilde A}$ is a ``line" for every ${\tilde A} \in {\tilde {\mathcal A}}$ which intersects ${\tilde S}$ at a ``plane". In $\Gamma$ this means that $x^{\perp} \cap A$ is a plane for every $A \in {\mathcal A}(S)$ which contains the point $p$.
\end{proof}

\begin{cor} Let $x \in {\mathcal X}(S)$ and $\lbrace p \rbrace = x^{\perp} \cap
S$. Let $y \in \nas{S}$ be such that $p \in y^{\perp} \cap S$. Then either $y \in x^{\perp}$ or $\lbrace x, y \rbrace$ is a polar pair.
\end{cor}

\begin{lem} Let $R, S \in {\mathcal S}$. If $R \cap S= \overline B$, a maximal
singular subspace in $\mss{\mathcal B}{S} \cap \mss{\mathcal B}{R}$ then $D(S) = D(R)$.
\end{lem}

\begin{proof} It suffices to prove $D(S) \subseteq D(R)$.

\medskip
  {\bf (i).} Let $x \in S$. If $x \in \overline B \subset R$ then $x \in D(R)$. So we may assume $x \in S \setminus \overline B$. Then $x ^{\perp} \cap \overline B$ is a plane and, by $(L.5)$, it follows that $x^{\perp} \cap R = \overline A_x$ is a maximal singular subspace of $R$ which belongs to $\mss{\mathcal A}{R}$ (since it has a plane in common with $\overline B$). Therefore $x \in \nas{R}$.

\medskip
 {\bf (ii).} Let now $x \in \nas{S}$. If $x \in R$ we are done. Let us assume that $x \not \in R$. Then $\overline A_x = x^{\perp} \cap S$ can have either a plane or a point in common with $\overline B$. If $\overline A_x \cap \overline B$ is a plane, then $x \in \nas{R}$ by a similar argument to that used in $(i)$.\\
Next let $\overline A_x \cap \overline B = \lbrace p \rbrace$, a single point.
First suppose that there exists a point $y \in x^{\perp} \cap R \setminus \lbrace p \rbrace$. {\it Claim}: $y \not \in A_x$. Assume to the contrary that $y \in A_x$. Therefore $y^{\perp} \cap S = {\overline A}_x$. Also, because $y \in R,\; y^{\perp} \cap S$ contains the plane $y^{\perp} \cap {\overline B}$. Then $y^{\perp} \cap S \supset \langle {\overline A}_x, y^{\perp} \cap {\overline B}\rangle $ and since ${\overline A}_x \cap {\overline B}$ is a single point, we get a contradiction with $(L.5)$. This proves the claim. Let $r \in {\overline A}_x \setminus y^{\perp}$, which, according to Step $(i)$, is in $\nas{R}$. Apply Lemma $3.1$ to the pair $\lbrace x, r \rbrace$ and symplecton $R$. It follows that $x \in \nas{R}$.\\
Next assume that $x^{\perp} \cap R = \lbrace p \rbrace$, a single point. {\it Claim}: $x \in {\mathcal X}(R)$. Let $q \in p^{\perp} \cap R$. Now $q \in R \subset \nas{S}$. Therefore the pair $\lbrace x, q \rbrace$ is polar.\\
For the remainder of this step we may assume that $q \not \in S$. Let ${\overline A}_q = q^{\perp} \cap S$ and $A_q$ be the maximal singular subspace containing ${\overline A}_q$. Since $x^{\perp} \cap A_q$ contains the line ${\overline A}_x \cap {\overline A}_q$, it follows, by $(L.4)$, that $x^{\perp} \cap A_q$ is a plane. This concludes the proof of the claim.

\medskip
{\bf (iii).} Let $x \in {\mathcal X}(S)$ be such that $x^{\perp} \cap {\overline B} = \lbrace p \rbrace$. Assume first $x^{\perp} \cap R = \lbrace p \rbrace$, a single point. Let $y \in p^{\perp} \cap R$. Since $y \in R$, as proved in Step $(i),\; y \in \nas{S}$ and, by Corollary $5.2, \; \lbrace x, y  \rbrace$ is a polar pair. Furthermore, if $A \in {\mathcal A}(S) \cap {\mathcal A}(R) \cap {\mathcal A}_p$ then $x^{\perp} \cap A$ is a plane. Therefore $x \in {\mathcal X}(R)$.\\
Assume now $x^{\perp} \cap R \supseteq py$, a line. Consider $\res{\Gamma}{p}$. We use the notations introduced in the proof of Lemma $5.1$. In $\res{\Gamma}{p}$, $\widetilde R$ and $\widetilde S$ are two ``symplecta" which meet at a plane of ${\mathcal B}$-type. Also $\tilde x$ is a ``point" at distance two from every single ``point" in $\widetilde S$ and $\tilde y \in \tilde x ^{\tilde \perp} \cap \widetilde R$. According to the result of Lemma $2.3$, it follows that $\tilde x^{\tilde \perp} \cap \widetilde R$ is a ``plane". Therefore, in $\Gamma$, $x^{\perp} \cap R$ is a maximal singular subspace of $R$. Now $(x^{\perp} \cap R) \cap \overline B = \lbrace p \rbrace$ which implies $x^{\perp} \cap R \in \mss{\mathcal A}{R}$. Thus $x \in \nas{R}$.

\medskip
{\bf (iv).} Let now $x \in {\mathcal X}(S)$ be such that $x^{\perp} \cap {\overline B} = \emptyset$ and let $\lbrace p \rbrace = x^{\perp} \cap S$. Recall that by Step $(i), \; p \in \nas{R}$. Let $A_p \in {\mathcal A}(R) \cap {\mathcal A}_p$ be the maximal singular subspace which contains $\langle p, p^{\perp} \cap R \rangle $. Note that $A_p \in {\mathcal A}(S)$ also, since $A_p \cap S = \langle p, p^{\perp} \cap {\overline B}\rangle $. By Lemma $5.1, \; x^{\perp} \cap A_p$ is a plane. Since $A_p$ has singular rank $5$ and $A_p \cap R$ has singular rank $3$, it follows that $x^{\perp} \cap A_p \cap R \not = \emptyset$. Note that $x^{\perp} \cap A_p \cap R$ cannot contain a line. If $x^{\perp} \cap A_p \cap R$ contains a line, since $A_p \cap R$ has singular rank $3$, this implies $x^{\perp} \cap {\overline B} \not= \emptyset$ which contradicts the hypothesis on $x$. Let $\lbrace q \rbrace = x^{\perp} \cap A_p \cap R$. {\it Claim}: $x \in {\mathcal X}(R)$. We first prove that $x^{\perp} \cap R = \lbrace q \rbrace$. Assume by contradiction that there is a point $r \in x^{\perp} \cap R \setminus \lbrace q \rbrace$. According to the above argument, $r \not \in A_p$. Then, by Lemma $3.1$, applied to the pair $\lbrace x, p \rbrace$ and symplecton $R$ we get $x \in \nas{R}$. But this implies $x^{\perp} \cap R \cap {\overline B} \not = \emptyset$, which contradicts the fact that $x^{\perp} \cap {\overline B} = \emptyset$. Therefore $x^{\perp} \cap R = \lbrace q \rbrace$. Since $x^{\perp} \cap A_p$ is a plane, it remains to prove that given any point $t \in q^{\perp} \cap R$, the pair $\lbrace x, t \rbrace$ is polar. This is clearly true for any $t \in q^{\perp} \cap {\overline B} = p^{\perp} \cap {\overline B}$; see the definition of ${\mathcal X}(S)$. So let us assume that $t \not \in {\overline B}$. Now $t \in \nas{S}$ and since $p^{\perp} \cap t^{\perp} \cap S$ contains a plane, the pair $\lbrace p, t \rbrace$ is polar. Set $T = \ll p, t \gg$ and note that $T \cap S = \langle p, p^{\perp} \cap t^{\perp} \cap S\rangle  \in \mss{\mathcal B}{S}$. Now $x \in {\mathcal X}(S)$, $x^{\perp} \cap S \cap T \not = \emptyset$ and $x^{\perp} \cap T \supset pq$ so, according to Step $(iii)$, $x \in \nas{T}$. Therefore $x^{\perp} \cap T \cap t^{\perp}$ contains a plane and this proves that $\lbrace x,t \rbrace$ is a polar pair.
\end{proof}

\begin{lem} Let $S, R \in {\mathcal S}$ be such that $S \cap R = L$, a line. Assume that $R
\cap {\mathcal X}(S) \not= \emptyset$. Then $D(S) = D(R)$.
\end{lem}

\begin{proof} Let $S, R \in {\mathcal S}$ be such that $S \cap R = L$, a line. Let
$x \in R \cap {\mathcal X}(S)$. Denote $x ^{\perp} \cap S = \lbrace p
\rbrace$ and note
that $ p \in L$. Let $A \in {\mathcal A}(S)$ be such that $L \subset A$. Then, according to $(L.2)$, $R \cap A$ can be a line or a maximal singular subspace of $R$. Since $x \in {\mathcal X}(S),\; x^{\perp} \cap A $ is a plane. Let $y \in x^{\perp} \cap A \setminus \lbrace p \rbrace$. Now $d(x, q)=2$ for any point $q \in L \setminus \lbrace p \rbrace$. Therefore $y \in R = \ll x, q \gg$. Then $A \cap R \supseteq \langle y, L \rangle $ and by $(L.2)$ it follows that $A \in {\mathcal A}(R)$. We proved that if $A \in {\mathcal A}(S)$ is such that $L \subset A$ then $A \in {\mathcal A}(R)$ as well.\\
Let now $A_1, A_2 \in {\mathcal A}(S) \cap {\mathcal A}(R)$ be two distinct maximal singular subspaces such that $L \subseteq A_1 \cap A_2$. The line $L$ is contained in the symplecton $S$ and so is the intersection of two maximal singular subspaces of $S$ of the same class (which can be either of the classes). Let $x_1 \in A_1 \cap S \setminus L$ and $x_2 \in A_2 \cap R \setminus L$. Set $T = \ll x_1, x_2 \gg$. Note that $T \cap S = \langle x_1, x_1^{\perp} \cap A_2 \cap S \rangle $ and $T \cap R = \langle x_2, x_2^{\perp} \cap A_1 \cap R \rangle $ both singular subspaces of ${\mathcal B}$-type. Then, by Lemma $5.3$: $D(S) = D(T) = D(R)$.
\end{proof}

\begin{cor} Let $x \in {\mathcal X}(S)$ and $\lbrace p \rbrace = x^{\perp} \cap S$. For $q \in p^{\perp} \cap S$ set $R = \ll x, q \gg$. Then $D(S) = D(R)$.
\end{cor}

\begin{proof} Let $S$ and $R$ be two symplecta as in the hypothesis. We claim that $R \cap S = pq$, a line. If $R \cap S$ contains a plane then $x^{\perp} \cap S$ contains a line, contradicting the properties of $x \in {\mathcal X}(S)$. Also $x \in {\mathcal X}(S) \cap R$. Now Lemma $5.4$ applies and the result follows.
\end{proof}

\begin{prop} Let $S \in {\mathcal S}$ then $D(S)$ is a subspace of $\Gamma$.
\end{prop}

\begin{proof}Let $x, y \in D(S)$ be two collinear points. We have to prove that $xy \subset D(S)$.

\medskip
  {\bf (i).} If $x, y \in S \cup \nas{S}$ the result follows from Lemma $3.2$.

\medskip
 {\bf (ii).} If $x \in {\mathcal X}(S)$ and $ y \in S$ then $\lbrace y \rbrace = x^{\perp} \cap S$. In this case $xy \subset S \cup {\mathcal X}(S)$ follows from the gamma space property of $\Gamma$.

\medskip
{\bf (iii).} Assume now that $x \in \nas{S}$ and $y \in {\mathcal X}(S)$.
Let ${\overline A}_x = x^{\perp} \cap S$ and $\lbrace p \rbrace = y^{\perp} \cap S$. Then $p \in {\overline A}_x$. If $p \not \in {\overline A}_x$, it follows, by Lemma $3.1$ that $y \in \nas{S}$, contrary to the assumption. Let $q \in p^{\perp} \cap S \setminus {\overline A}_x$. Set $R = \ll q, x \gg$. Then, by Lemma $5.3$, $D(S) = D(R)$. Now $y \in {\mathcal X}(S) \subset D(R)$ and $y^{\perp} \cap R \supseteq xp$ implies $y \in \nas{R}$. Note that $y \not \in R$ because this would imply $y^{\perp} \cap S$ contains a plane. So $x \in R,\; y \in \nas{R}$ and, according to Lemma $3.2, \; xy \subset D(R) = D(S)$.

\medskip
{\bf (iv).} Let $x, y \in {\mathcal X}(S)$. First assume $x^{\perp} \cap
y^{\perp} \cap S = \lbrace p \rbrace$. Let $q \in p^{\perp} \cap S$ and set $R = \ll x, q \gg$. Then, according to Corollary $5.5$, $D(S) = D(R)$. Thus $x \in R$ and $y \in R \cup \nas{R}$ implies $xy \subset D(R)= D(S)$.\\
Let now $x^{\perp} \cap S = \lbrace p \rbrace$ and $y^{\perp} \cap R = \lbrace q \rbrace$ be such that $p \not = q$. If $p \in q^{\perp}$ then set $R = \ll x, q \gg$ and, by Corollary $5.5$, $D(R) = D(S)$. Therefore $xy \subset R \subset D(S)$. If $p \not \in q^{\perp}$ then take $r \in p^{\perp} \cap q^{\perp}$ and set $T = \ll x, r \gg$. Again, by Corollary $5.5,\; D(T) = D(S)$. So $x \in T,\; y \in D(T)$ and by the previous results of this Proposition, $xy \subset D(T) = D(S)$.
\end{proof}

In the next proof we use the assumption that $\Gamma$ satisfies the Weak Hexagon Axiom ({\it WHA}).
\begin{prop} Let $S \in {\mathcal S}$, then $D(S)$ is a $2$-convex
subspace of $\Gamma$.
\end{prop}

\begin{proof}Let $x, y \in D(S)$ be two points at distance $2$. We have to prove that $x^{\perp} \cap y^{\perp} \subset D(S)$.

\medskip
  {\bf (i).} If $x, y \in S$ then obviously $x^{\perp} \cap y^{\perp} \subset
S$ since $S$ is a $2$-convex subspace of $\Gamma$.

\medskip
 {\bf (ii).} Let now $x \in \nas{S}$ and $y \in S$. The pair $\lbrace x, y \rbrace$ is polar. Set $R = \ll x, y \gg$. Since $R \cap S \in \mss{\mathcal B}{S}$ it follows, by Lemma $5.3$ that $D(S) = D(R)$. Therefore $x^{\perp} \cap y^{\perp} \subset R \subset D(S)$.

\medskip
{\bf (iii).} Consider now the case when $x, y \in \nas{S}$. Then ${\overline A}_x = x^{\perp} \cap S$ and ${\overline A}_y =
y^{\perp} \cap S$ are two maximal singular subspaces of $S$ from the same
class. There are two cases to consider:\\
{\bf (iii.a).} Assume ${\overline A}_x \cap {\overline A}_y =L$, a line. Let $p \in {\overline A}_y \setminus L$ and set $R = \ll x,p \gg$. Then $R \cap S = \langle p, p^{\perp} \cap {\overline A}_x \rangle  \in M_{\mathcal B}(S)$ and, by Lemma $5.3,\; D(S) = D(R)$. So $y \in D(R)$ and because $y^{\perp} \cap R \supseteq \langle p, L \rangle $ it follows $y \in \nas{R}$. Note that $y \not \in R$ since $A_y \cap R = \langle p, x^{\perp} \cap A_y\rangle $ has singular rank $3$ and thus is already maximal in $R$. Now $x \in R,\; y \in \nas{R}$ and according to Step $(ii)$ of this Proposition, $x^{\perp} \cap y^{\perp} \subset D(R) =D(S)$.

\medskip
{\bf (iii.b).} Assume now that ${\overline A}_x$ and ${\overline A}_y$ are disjoint maximal singular subspaces of $S$. Let $z \in x^{\perp} \cap y^{\perp}$. Clearly $z \not \in S$. If $z^{\perp} \cap S \not = \emptyset$ then $z \in \nas{S}$, by Lemma $3.1$, and we are done. So we have to prove that $z^{\perp} \cap S \not = \emptyset$. Assume by contradiction that $z^{\perp} \cap S = \emptyset$. Note that this also implies $z^{\perp} \cap A_x = \lbrace x \rbrace$ and $z^{\perp} \cap A_y = \lbrace y \rbrace$, see $(L.4)$ and recall that ${\mathcal A}$-spaces have singular rank $5$. Take $x_1 \in {\overline A}_x$ and $y_1 \in {\overline A}_y\setminus x_1^{\perp}$. Also let $z_1 \in x_1^{\perp} \cap y_1^{\perp}$ be such that $z_1 \not \in x^{\perp} \cup y^{\perp}$. Now $H = (x_1, z_1, y_1, y,z,x)$ is a minimal $6$-circuit in $\Gamma$ which contains at least one polar pair $\lbrace x_1, y_1 \rbrace$. Then ({\it WHA}) applies and there exists a point $w \in z^{\perp} \cap x_1^{\perp} \cap y_1^{\perp} \subset S$. But this contradicts the assumption made and it follows that $z^{\perp} \cap S \not = \emptyset$.

\medskip
{\bf (iv).} Let $x \in {\mathcal X}(S)$ and $y \in S$. Let $\lbrace p \rbrace = x^{\perp} \cap S$. If $y \in p^{\perp}$ then set $R = \ll x,y \gg$ which, by Corollary $5.5$ is such that $D(R) = D(S)$ and $x^{\perp} \cap y^{\perp} \subset R \subset D(S)$. So let us assume $y \not \in p^{\perp}$. Take $A \in {\mathcal A}(S) \cap {\mathcal A}_p$. Then according to Lemma $5.1$, $x^{\perp} \cap A$ is a plane. Let $x_1 \in x^{\perp} \cap A \setminus \lbrace p \rbrace$. Also let $r \in p^{\perp} \cap S \setminus A$. Set $R = \ll x_1, r \gg$. Since $R \cap S = \langle r, r^{\perp} \cap A \rangle  \in \mss{\mathcal B}{R} \cap \mss{\mathcal B}{S}$ it follows, by Lemma $5.3$ that $D(S) = D(R)$. Now $x^{\perp} \cap R \supset x_1p$ so $x \in \nas{R}$. Also $y \in \nas{R}$. According to Step $(iii)$ of this Proposition $x^{\perp} \cap y^{\perp} \subset D(R) = D(S)$.

\medskip
{\bf (v).} Let $x \in \nas{S}, \; y \in {\mathcal X}(S)$. Let
$\lbrace p \rbrace = y^{\perp} \cap S$ and ${\overline A}_x = x^{\perp} \cap S$. Assume first that $p \not\in {\overline A}_x$. Set $R = \ll x, p \gg$. Then $R \cap S
\in \mss{\mathcal B}{S}$ and according to Lemma $5.3, \; D(R) = D(S)$. Now $x \in R,\; y \in D(R)$ and using the above results of this Proposition, it follows $x^{\perp} \cap y^{\perp} \subset D(R) = D(S)$.\\
Next consider the case $p \in {\overline A}_x$. Let $q \in {\overline A}_x$ be a point and set $T = \ll q,y \gg$. By Corollary $5.5,\; D(S) = D(T)$. So we may apply the results of the Steps $(i)-(iii)$ of this Proposition to the pair of points $ x\in D(T),\; y \in T$ and conclude that $x^{\perp} \cap y^{\perp} \subset D(T) = D(S)$.

\medskip
{\bf (vi).} Let now $x, y \in {\mathcal X}(S)$ be such that $x^{\perp} \cap S =
\lbrace p \rbrace$ and $y^{\perp} \cap S = \lbrace q \rbrace$. If $p \not \in q^{\perp}$ then let $r \in p^{\perp} \cap q^{\perp}$; if $p \in q^{\perp} \setminus \lbrace q \rbrace$ then we let $r = q$ and if $p = q$ then take $r \in p^{\perp}$. Set $R =  \ll x, r \gg$. By Corollary $5.5,\; D(S) = D(R)$. Now $y \in {\mathcal X}(S) \subset D(R), \; x \in R$ and using the previous results of this Proposition, $x^{\perp} \cap y ^{\perp} \subset D(R) = D(S)$.
\end{proof}

\begin{lem} Let $x \in {\mathcal X}(S)$ with $\lbrace p \rbrace = x^{\perp} \cap S$ and let $z \in S$. Then $d(x, z) = 1+d(p, z)$.
\end{lem}

\begin{proof} We start by proving the following {\it claim}: if $x, y \in {\mathcal X}(S)$ are two collinear points with $\lbrace p \rbrace = x^{\perp} \cap S$ and $\lbrace q \rbrace = y^{\perp} \cap S$ then $q \in p^{\perp}$. Assume by contradiction that $q \not \in p^{\perp}$. Let $A \in {\mathcal A}_p \cap{\mathcal A}(S)$. Then, according to Lemma $5.1,\; x^{\perp} \cap A$ is a plane. Let $x_1 \in x^{\perp} \cap A \setminus \lbrace p \rbrace$. Also let $r \in p^{\perp} \cap q^{\perp} \setminus A$. Set $R = \ll r, x_1 \gg$. Now $R \cap S = \langle r, r^{\perp} \cap A \cap S\rangle  \in \mss{\mathcal B}{S}$ and, by Lemma $5.3, \; D(R) = D(S)$. Note that $x^{\perp} \cap R \supset px_1$ and consequently $x \in \nas{R}$. Set $\overline A_x = x^{\perp} \cap R$. Also $y \in {\mathcal X}(S) \subset D(R)$ implies that there is a point $y_1 \in y^{\perp} \cap R$. If $y_1 \not \in {\overline A}_x$, since $x$ and $y$ are collinear, it follows, by Lemma $3.1$, that $y \in \nas{R}$. But then $y^{\perp} \cap R$ and $R \cap S$ are maximal singular subspaces in $R$ from different families and therefore they have a point in common. Since $q \not \in R \cap S$ it follows that $y^{\perp} \cap R$ contains more than a point, which is a contradiction with the fact that $y \in {\mathcal X}(S)$. Therefore $y_1 \in {\overline A}_x$. Now $q \in S$ and, by Lemma $5.3$, $q \in \nas{R}$. Let $\overline A _q = q^{\perp} \cap R$. But $q \in y^{\perp}$ so $y_1 \in {\overline A}_q$ because otherwise another application of Lemma $3.1$ to the pair $\lbrace y, q \rbrace$ and symplecton $R$ gives $y \in \nas{R}$, a contradiction. Therefore $y_1 \in {\overline A}_x \cap {\overline A}_q$ which, using $(L.3)$ implies ${\overline A}_x \cap {\overline A}_q = L$ a line. But $q \not \in p^{\perp}$ so $p \not \in L$. Since ${\overline A}_q$ has singular rank $3$ and since $S \cap R$ meets ${\overline A}_q$ in a plane it follows that $L \cap (R \cap S) \not = \emptyset$. But this implies that $x^{\perp} \cap S$ contains more than a point. We reach a contradiction with the fact that $x \in {\mathcal X}(S)$. Therefore the assumption made was false and $q \in p^{\perp}$.

\medskip
In order to prove the Lemma it suffices to prove that, if $x \in {\mathcal X}(S)$ and $z \in S \setminus p^{\perp}$ then $d(x, z) = 3$. Assume by contradiction that $d(x, z) = 2$. Then there exists a point $t \in x^{\perp} \cap z^{\perp}$. Note that $t$ cannot be in $S$. According to Proposition $5.7,\; t \in D(S)$. Moreover $t \not \in \nas{S}$, because since $z \not \in p^{\perp}$, Lemma $3.1$ would imply $x \in \nas{S}$, a contradiction. So we must have $t \in {\mathcal X}(S)$. But now, according to the previous paragraph $z \in p^{\perp} \cap S$, which contradicts the hypothesis on $z$. Therefore the assumption made was false and in this case $d(x,z)=3$.
\end{proof}

\begin{prop} For any $S \in {\mathcal S},\ D(S)$ is a strong parapolar
subspace of $\Gamma$.
\end{prop}

\begin{proof} We have to prove that if $x,y \in D(S)$ are two points at distance two, the pair $\lbrace x,y \rbrace$ is polar, that is $x^{\perp} \cap y^{\perp}$ contains at least two points.

\medskip
  {\bf (i).} Let $x \in D(S)$ and $y \in S$. If both $x$ and $y$ are in $S$ then $S = \ll x, y \gg$ and we are done. Assume next $x \in \nas{S}$. Then $y^{\perp} \cap x^{\perp} \cap S$ contains a plane and therefore $\lbrace x, y \rbrace$ is a polar pair. Let now $x \in {\mathcal X}(S)$ with $\lbrace p \rbrace = x^{\perp} \cap S$. In this case, since $d(x,y) = 2$, by Lemma $5.8$ it follows that $y \in p^{\perp}$. Then the fact that $\lbrace x, y \rbrace$ is polar pair follows from the definition of ${\mathcal X}(S)$.

\medskip
{\bf (ii).} Assume now $x, y \in \nas{S}$. Let ${\overline A}_x = x^{\perp} \cap S$ and ${\overline A}_y = y^{\perp} \cap S$. If ${\overline A}_x \cap {\overline A}_y \not = \emptyset$ the result is immediate. So we may assume that ${\overline A}_x \cap {\overline A}_y = \emptyset$. Let $p \in {\overline A}_x$. Then set $R = \ll p, y \gg$ which, by Lemma $5.3$, is such that $D(R) = D(S)$. Thus $x \in {\mathcal X}(S) \subset D(R),\; y \in R$ and by Step $(i)$ of this Proposition it follows that $\lbrace x, y \rbrace$ is a polar pair.

\medskip
{\bf (iii).} Let $x \in \nas{S}$ and $y \in {\mathcal X}(S)$. Let ${\overline A}_x = x^{\perp} \cap S$ and $\lbrace p \rbrace = y^{\perp} \cap S$. If $p \in {\overline A}_x$ then, according to Corollary $5.2,\; \lbrace x, y \rbrace$ is a polar pair. Assume next that $p \not \in {\overline A}_x$. Take a point $q \in p^{\perp} \cap {\overline A}_x$ and set $R = \ll y, q \gg$. By Corollary $5.5, \; D(R) = D(S)$. So $y \in R,\; x \in \nas{S} \subset D(R)$ and $\lbrace x, y \rbrace$ is a polar pair by Step $(i)$ of this Proposition.

\medskip
{\bf (iv).} Let now $x, y \in {\mathcal X}(S)$ with $\lbrace p \rbrace = x^{\perp} \cap S$ and $\lbrace q \rbrace = y^{\perp} \cap S$. If $p \in q^{\perp} \setminus \lbrace q \rbrace$ set $R = \ll p, y \gg$. If $p = q$ then take $r \in p^{\perp}$ and if $p \not \in q^{\perp}$ take $r \in p^{\perp} \cap q^{\perp}$ and set $R =\ll r, y \gg$. Then, by Corollary $5.5,\; D(R) = D(S)$ and since $x \in {\mathcal X}(S) \subset D(R), y \in R$ then $\lbrace x, y \rbrace$ is a polar pair, see Step $(i)$ of this Proposition again.
\end{proof}
\bigskip

Let $S \in {\mathcal S}$. Then $D(S) \in {\mathcal D}$ is a $2$-convex
strong parapolar subspace of $\Gamma$, that is if $x, y \in D(S)$ are two points at distance two, then $R = \ll x, y \gg \subset D(S)$. Next, we analyze the possible relations between $D(S)$ and $D(R)$.\\

{\bf 1.} Let $x \in D(S)$ and $y \in S$. If $x \in S$ then $S= \ll x,y \gg = R$. Next let $x \in \nas{S}$ and set ${\overline A}_x = x^{\perp} \cap S$. Then $R \cap S = \langle y, y^{\perp} \cap {\overline A}_x \rangle  \in \mss{\mathcal B}{S}$. Therefore $D(S) = D(R)$, by Lemma $5.3$. Let now $x \in {\mathcal X}(S)$ with $\lbrace p  \rbrace = x^{\perp} \cap S$. By Lemma $5.8$, $y \in p^{\perp} \cap S$ and by Corollary $5.5$ it follows $D(R) = D(S)$.

\medskip
{\bf 2.} Let $x, y \in \nas{S}$ with ${\overline A}_x = x^{\perp} \cap S$ and ${\overline A}_y = y^{\perp} \cap S$.\\
{\bf 2.a.} Assume ${\overline A}_x \cap {\overline A}_y = L$, a line. {\it Claim}: $R \cap S = L$. Suppose by contradiction that $R \cap S$ contains a plane $\langle r,L \rangle $ where $r$ is a point not on $L$. Since $r \not \in A_x \cap A_y$ we can assume, without loss of generality, that $r \not \in A_y$. If $r \in {\overline A}_x$ then $\langle r, y^{\perp} \cap A_x, x \rangle$ is a singular subspace of rank $4$ in $R$, which has polar rank $4$, a contradiction. So let us assume that $r \in S \setminus ({\overline A}_x \cup {\overline A}_y)$. Then $\langle r, r^{\perp} \cap {\overline A}_x, r^{\perp} \cap {\overline A}_y\rangle  \subset R \cap S$ which is a contradiction with the fact that $R \cap S$ can be at most a common maximal singular subspace. Therefore the claim is proved.\\
Now pick a point $z \in x^{\perp} \cap y^{\perp} \setminus L$ such that $z^{\perp} \cap L = \lbrace p \rbrace$, a single point. We intend to prove that $z \in {\mathcal X}(S)$. According to Proposition $5.7, \; z\in D(S)$, so it suffices to prove that $z^{\perp} \cap S = \lbrace p \rbrace$. Assume by contradiction that $z \in \nas{S}$ with ${\overline A}_z = z^{\perp} \cap S$. Then, according to $(L3)$, ${\overline A}_z \cap {\overline A}_x$ is a line distinct from $L$. So if $t \in L \setminus \lbrace p \rbrace$ then we may write $R = \ll z, t \gg$. Then ${\overline A}_z \cap {\overline A}_x \subset z^{\perp} \cap t^{\perp} \subset R$ and we reach a contradiction with the previous result that $R \cap S = L$. Therefore $z^{\perp} \cap S = \lbrace p \rbrace$ and $z \in {\mathcal X}(S)$. Now $R$ and $S$ are as in Lemma $5.4$ and consequently $D(R) = D(S)$.\\
{\bf 2.b.} Let ${\overline A}_x \cap {\overline A}_y = \emptyset$. Let $t \in {\overline A}_x$ and set $T = \ll t,y \gg$. Thus $T \cap S = \langle  t, t^{\perp} \cap {\overline A}_y \rangle  \in \mss{\mathcal B}{S}$ and by Lemma $5.3,\; D(T) = D(S)$. Now $x \in D(T), y \in T$ so by Step $1$, $R = \ll x, y\gg$ is such that $D(R) = D(T) = D(S)$.\\

{\bf 3.} Let $x \in \nas{S}, \; y \in {\mathcal X}(S)$ with $y^{\perp} \cap S = \lbrace p \rbrace$. If $p \in {\overline A}_x$ then let another point $q \in {\overline A}_x$ and set $T = \ll q, y \gg$. By Corollary $5.5,\; D(T) = D(S)$. Apply Step $1$ to $x \in D(T),\; y \in T$ and get $D(R) = D(T) = D(S)$. If $p \not \in {\overline A}_x$ set $T = \ll p, x \gg$ and since $T \cap S \in \mss{\mathcal B}{S}$, Lemma $5.3$ gives $D(T)=D(S)$. Now $x \in T$ and $y \in D(T)$ and using Step $1$ one more time we get $D(R) = D(T) = D(S)$.\\

{\bf 4.} Let $x, y \in {\mathcal X}(S)$ with $\lbrace p \rbrace =
x^{\perp} \cap S$ and $\lbrace q \rbrace = y^{\perp} \cap S$. If $p = q$ take $r \in p^{\perp} \cap S$, set $T = \ll r, x \gg$ and apply Step $1$. If $p \in q^{\perp} \setminus \lbrace q \rbrace$ then take $T = \ll x, q \gg$ and apply Step $1$. If $p \not \in q^{\perp}$ then take $r \in p^{\perp} \cap q^{\perp}$, set $T =  \ll x,
r \gg$ and apply Step $1$ again. Then $D(R) = D(T) = D(S)$.\\

Therefore, given $S \in {\mathcal S},\ D(S) \in {\mathcal D}$ and some
symplecton $R \subset D(S)$, we see that $D(R) = D(S)$ and we have proved the
following:
\begin{prop} Given $R \in {\mathcal S}$ there exists a unique element
$D(S) \in {\mathcal D}$, for some $S \in {\mathcal S}$, containing $R$.
\end{prop}
\begin{prop} Let $D(S) \in {\mathcal D}$ for some $S \in {\mathcal S}$ and
let $(x, R)$ a non-incident point-symplecton pair in $D(S)$. Then
$x^{\perp} \cap R$ is a point or a maximal singular subspace of $R$.
\end{prop}
\begin{proof} Since $R \subset D(S)$, by the previous Proposition $5.10, \; D(R) = D(S)$. If $x \not \in R$ then $x \in \nas{R} \cap
{\mathcal X}(R)$ and the conclusion follows at once.
\end{proof}
\bigskip
\begin{proof}[{\bf Proof of the Theorem 2}] Let {\pl} be a parapolar space
which is locally $A_{n-1,4}({\mathbb K})$ for $n =7$ or $n$ an integer greater than $8$. Assume that $\Gamma$ satisfies the Weak Hexagon Axiom ({\it WHA}). Given $S \in {\mathcal S}$ a symplecton, we set $D(S) = S \cup \nas{S} \cup {\mathcal X}(S)$ with $\nas{S}$ and ${\mathcal X}(S)$ defined as in Section $3$. Then $D(S)$ is a $2$-convex subspace of
$\Gamma$, see Propositions $5.6$ and $5.7$, which is strong parapolar, by
Proposition $5.9$, and such that, for any non-incident point-symplecton
pair $(x, R)$ in $D(S)$, $x^{\perp} \cap R$ is a point or a maximal
singular subspace of $R$. By the characterization theorem of Cohen and
Cooperstein \cite{cc}, see section $2.6$ also, $D(S)$ is isomorphic to the half-spin geometry $D_{6,6}({\mathbb K})$.
Moreover, by Proposition $5.10$, every symplecton lies in a unique element
of ${\mathcal D}$. This ends the proof of the Theorem.
\end{proof}

\section{The sheaf construction}
Let $n$ be an integer greater than $6$. Let $I = \lbrace
1, \dots,n \rbrace$. For $1 <  k \leq n-4$, let $K = \lbrace k, \dots, k+4
\rbrace$ and set $J = I \setminus K$, $J_1 = \lbrace 1, \dots, k-1 \rbrace,\;  J_2 =
\lbrace k+5, \dots, n \rbrace$. Let $D$ be the following locally truncated diagram with typeset $I$:\\
{\it 
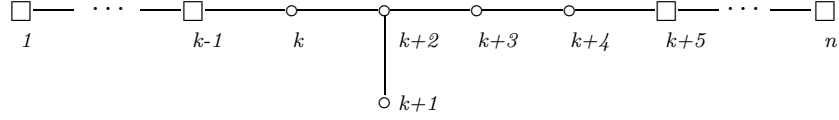
\begin{picture}(1000, 50)(0,0)
\put(50,25){$\qed$}
\put(70,29){\line(1,0){15}  \; $\dots$ \line(1,0){15}}
\put(65,15){\scriptsize 1}
\put(115,25){$\qed$}
\put(130,15){\scriptsize k-1}
\put(135,29){\line(1,0){30}}
\put(165,26){$\circ$}
\put(168,15){\scriptsize k}
\put(170,29){\line(1,0){30}}
\put(200,26){$\circ$}
\put(208,15){\scriptsize k+2}
\put(203,-3){\line(0,1){30}}
\put(200,-9){$\circ$}
\put(208,-9){\scriptsize k+1}
\put(205,29){\line(1,0){30}}
\put(235,26){$\circ$}
\put(238,15){\scriptsize k+3}
\put(240,29){\line(1,0){30}}
\put(270,26){$\circ$}
\put(273,15){\scriptsize k+4}
\put(275,29){\line(1,0){30}}
\put(294,25){$\qed$}
\put(309,15){\scriptsize k+5}
\put(314,29){\line(1,0){15} \;$\ldots$ \line(1,0){15}}
\put(354,25){$\qed$}
\put(369,15){\scriptsize n}
\end{picture}}
\vspace{.5cm}

Suppose $\Gamma$ is a geometry over $K$ which is $J$-locally truncated
with respect to the diagram $D$ over $I$ and let ${\mathcal F}$ be the set
of all nonempty flags of $\Gamma$. The existence of a sheaf for such a locally
truncated geometry was stated without proof in Brouwer and Cohen (Theorem $4$, \cite{loc}). In this Section we give a detailed construction of this sheaf. See Subsection $2.4$ for the definition of the sheaf. Before proceeding with the sheaf construction, we need one more result. For completeness we also give its proof.
\newtheorem*{ela}{Lemma}
\begin{ela}(Ellard and Shult \cite{els}) Let $D$ be a diagram over $I$. Let
$K$ be a subset of $I$, $|K|>3$, and suppose
that $\Gamma$ is a locally truncated geometry over $K$ with diagram
$D$. Assume that $\Sigma$ is a sheaf over
${\mathcal F}_2$, the family of all non-empty flags of $\Gamma$ of rank at most two having typeset $I$
and such that for each $F \in {\mathcal F}_2$, $\res{\Gamma}{F} = \; _J \Sigma
(F)$ where $J = I \setminus K$. Then $\Sigma$ can be extended to a sheaf $\Sigma '$ with typeset $I$ defined over the family ${\mathcal
F}$ of all non-empty flags of $\Gamma$, enjoying the same extended identity on
$J$-truncations: $_J \Sigma '(F) = \res{\Gamma} {F}$, for all $F \in {\mathcal F}$.
\end{ela}
\begin{proof} Set $\Sigma'(x) = \Sigma(x)$ for each object $x$ in $\Gamma$. Since $\Sigma$ is a sheaf , we have:
\begin{center}
$\Sigma ' (x, y) = \res{\Sigma (x)}{y} = \res{\Sigma (y)}{x}$
\end{center}
for all incident pairs of distinct objects $\lbrace x, y \rbrace$ of $\Gamma$. This defines $\Sigma' (F)$ for all flags $F = \lbrace x, y \rbrace$ of rank $2$.\\
Now let $F$ be any flag of rank greater than $2$ and define:
\begin{center}
$\Sigma '(F) := \res{\Sigma (x)}{F \setminus \lbrace x \rbrace}$
\end{center}
where $x$ is some object in $F$. To show that $\Sigma '(F)$ is independent on the choice of $x$, suppose $y$ is a second object in $F$. Then:
\begin{equation*}
\begin{split}
\Sigma ' (F) & = \res{\Sigma(x)}{F \setminus \lbrace x \rbrace} = \res{ \res{\Sigma (x)}{y}}{ F \setminus \lbrace x, y \rbrace } =\\
& = \res{ \res{\Sigma (y)}{x}}{ F \setminus \lbrace x, y \rbrace } = \res{\Sigma(y)}{\lbrace x \rbrace \cup (F \setminus \lbrace x, y \rbrace )}=\\
& = \res{\Sigma (y)}{F \setminus \lbrace y \rbrace}
\end{split}
\end{equation*}
Therefore $\Sigma'(F)$ is well defined. Clearly, for any such a flag
\begin{center}
$_J \Sigma ' (F) =\;  _J \res{\Sigma '(x)}{F \setminus \lbrace x \rbrace} = \res {\res {\Gamma}{x}}{F \setminus \lbrace x \rbrace} = \res {\Gamma} {F}$
\end{center}
since the processes of truncation and taking residuals commute.\\
Next show that $\Sigma'$ is a sheaf. Let $F$ and $F_1$ be flags of $\Gamma$ with $\emptyset \not= F \subseteq F_1$. We must show that
\begin{center}
$\Sigma' (F_1) = \res {\Sigma' (F)}{F_1 \setminus F}$
\end{center}
Choose $x \in F$. Then $x \in F_1$, so by definition
\begin{center}
$\Sigma'(F_1) = \res {\res {\Sigma (x)}{F \setminus \lbrace x \rbrace}}{F_1 \setminus F} = \res {\Sigma' (F)}{F_1 \setminus F}$
\end{center}
Finally, let $F_1 \subseteq F_2 \subseteq F_3$, then:
\begin{equation*}
\begin{split}
\Sigma'(F_3) & = \res{\Sigma'(F_2)}{F_3 \setminus F_2}=\res{\res{\Sigma'(F_1)}{F_2 \setminus F_1}}{F_3 \setminus F_2}=\\
& = \res{\Sigma'(F_1)}{(F_2 \setminus F_1) \cup (F_3 \setminus F_2)} = \res{\Sigma'(F_1)}{F_3 \setminus F_1}
\end{split}
\end{equation*}
The proof is complete.\\
\end{proof}

Therefore, it will suffice to construct the sheaf over the family of flags
of rank $1$ and $2$ of $\Gamma$. In what follows we give the step by step construction of the sheaf:\\

{\bf 1}. Let $a \in {}^{k}\Gamma$. Recall that this means $a$ is an object of type $k$ in $\Gamma$. We define $\Sigma ^R(a)$ to be a
geometry
of locally truncated type belonging to the diagram $D_{I \setminus
(\lbrace k \rbrace
\cup J_1)}$ which satisfies the property that $\res {\Gamma} {a} =\;  _{J_2}\Sigma ^R
(a)$. Note
that this is well defined since $\Sigma ^R (a)$ is of $J_2$-locally type $A_{n-k, 2}$,
which, up to the relabeling of the nodes is unique (see Brouwer and Cohen \cite{loc}). In $\Sigma
^R (a)$ the
objects of types in $K$ are the same as in $\Gamma$ with their corresponding
incidence. The objects of type $i \in J_2$ are collections of objects in $\Gamma$
with their flags, which are incident with a given object of type in $K$; the incidence between objects of types in $J_2$ is given by the symmetrized containment.\\

{\bf 2}. Let $b \in \; ^{k+1}\Gamma$. Then $\res{\Gamma}{b}$ is
the $J$-truncation of a geometry belonging to the diagram $D_{I
\setminus \lbrace k+1 \rbrace}$ of type $A_{n-1, k+1}$. But this is
uniquely determined by its truncation and thus we may unambiguously define
$\Sigma(b)$
to be such that $\res{\Gamma}{ b} =\;  _J \Sigma(b)$. There are two types of objects in $\Sigma
(b)$: those inherited from $\Gamma$, with their incidence and the objects with types
in $J \setminus K$, which are defined as in Step $1$.\\

{\bf 3}. Let now $l \in \lbrace k+1, \dots, k+4 \rbrace$. Let $x_l \in \;
^{l}\Gamma$ and $F = \lbrace a,x \rbrace $ a $\lbrace k, l-1 \rbrace $ - flag in $\res{\Gamma}{x_l}$. When $l=k+1$, $F = \lbrace a \rbrace$ is just an object of type $k$. Define
recursively:
\begin{center}
$\Sigma ^R (x_l) = \Sigma ^R (a, x_l) := \res{\Sigma ^R (a,x)}{x_l}$
\end{center}
In order to prove that this is well defined, consider $F'$ another flag of type
$\lbrace k, l-1 \rbrace$ in $\res{\Gamma }{x_l}$, which is $i$-adjacent to $F$ in $C$, the chamber system associated to $^{\lbrace k, l-1 \rbrace}\res {\Gamma} {x_l}$ with  $i \in \lbrace k, l-1 \rbrace$. Then:
\begin{center}
$\res{\Sigma ^R (F)}{x_l} = \res{\Sigma ^R (F \cap F')} {(F \setminus F') \cup
\lbrace x_l \rbrace } = \res{\Sigma ^R (F') }{x_l}$
\end{center}
which proves the well-definedness in this case. If $F$ and $F'$ are not
$i$-adjacent, since the chamber system $C$ is connected, there is a
gallery from $F$ to $F'$, and by repeated applications of the above
argument we get the result.\\

{\bf 4}. Let $e \in \; ^{k+4} \Gamma$. Define $\Sigma ^L (e)$ to be the
geometry of $J_1$-locally truncated type, with diagram $D_{I \setminus (\lbrace
k+4 \rbrace \cup J_2 )}$, and such that $\res{\Gamma}{ e} = \;
_{J_1} \Sigma ^L (e)$. This is uniquely defined, since
it has truncated type $D_{k+2, k+2}$. The objects and their incidence in
$\Sigma^L(e)$ are defined in similar manner to those for $\Sigma ^R (a)$,
Step $1$.\\

{\bf 5}. Let $l \in \lbrace k+3, \dots, k \}$ (taken in this order) and define the $L$ (left)
counterparts of the quantities introduced in Step $3$. Let $x_l \in
\; ^{l}\Gamma$ and $F = \lbrace x, e \rbrace$ a $\lbrace l+1, k+4 \rbrace$-flag in $\res{\Gamma}
{x_l}$. Recursively define:
\begin{center}
$\Sigma ^L (x_l) = \Sigma ^L (x_l, e) := \res{\Sigma ^L(x, e)} {x_l}$
\end{center}
which can be proved that it is a good definition in the same way as in
Step $3$.\\

The sheaf values over the rank $1$ flags of $\Gamma$ can be written as
follows: \begin{itemize}
\item[{\bf (i).}] \quad $\Sigma (x) := \Sigma ^L (x) \oplus \Sigma ^R (x),\;   \text{for
any object}\ x\  \text{of type in}\  K \setminus \lbrace k+1 \rbrace$;
\item[{\bf (ii).}] \quad $\Sigma (b) \; \text{for any object}\ b \in \; ^{k+1}\Gamma$.
\end{itemize}

Note that in the case of the objects $c \in \; ^{k+2}\Gamma$, the definition gives:
\begin{center}
$\Sigma(c) = \Sigma ^L (b) \oplus \; ^{k+1}\Gamma \oplus \Sigma ^R (c)$
\end{center}
where $b$ is an object of type $k+1$, incident with $c$.\\

We proceed now with the construction of the sheaf over the rank $2$
flags of $\Gamma$. In this case, the sheaf values should agree with the
values calculated at the corresponding residuals, this means that, if
$ \lbrace x,y\rbrace $ is a nonempty rank $2$ flag of $\Gamma$, then the following must
be true: \begin{itemize}
\item[{\bf (S).}] \quad $\res{\Sigma(x)}{y} = \Sigma(x,y) = \res{\Sigma(y)}{x}$
\end{itemize}

{\bf 6}. Let $\lbrace x_i, x_j \rbrace $ be a rank $2$ flag of type $\lbrace i, j
\rbrace$ in $\Gamma$ with $i, j \in K \setminus \lbrace k+1 \rbrace$
and $i<j$. The corresponding sheaf value is defined to be:
\begin{center}
$ \Sigma (x_i, x_j) := \Sigma ^{L} (x_i) \oplus \;
^{(ij)} \Gamma \oplus \Sigma ^{R}(x_j)$
\end{center}
where $^{(ij)}\Gamma$ denotes the
truncation of $\Gamma$ to those objects of type $l \in K$ with $i<l<j$. If
$i = j-1$ then $^{(ij)}\Gamma$ is empty. Next we check the property $(S)$:
\begin{equation*}
\begin{split}
\res{\Sigma(x_i)}{x_j} &= \res{\Sigma^L (x_i) \oplus \Sigma ^R (x_i)}{x_j} = \res {\Sigma ^L (x_i)}{x_j} \oplus \res {\Sigma ^R(x_i)}{x_j} =\\
& = \Sigma ^L (x_i) \oplus \; ^{(ij)}\Gamma \oplus \Sigma ^R (x_j) = \res {\Sigma ^L (x_j)}{x_i} \oplus \Sigma ^R (x_j) = \\
& = \res {\Sigma ^L (x_j) \oplus \Sigma ^R (x_j)} {x_i} = \res{\Sigma (x_j)}{x_i}
\end{split}
\end{equation*}

{\bf 7}. Consider now those rank $2$ flags $\lbrace x_i, x_j \rbrace $ in $\Gamma$, in
which one of the types $i$ or $j$ is $k+1$, the other type taking any
other value in $K$.\\
 {\bf 7.i.} \quad Let $\lbrace a,b \rbrace $ be a $\{ k, k+1 \}$ flag. Then define:
\begin{center}
$\Sigma (a,b) := \Sigma ^L(a) \oplus \Sigma ^R (b)$
\end{center}
and check the sheaf property $(S)$:
\begin{equation*}
\begin{split}
\res{\Sigma(a)}{b} & = \res{\Sigma ^L(a) \oplus \Sigma ^R (a)}{b} = \res {\Sigma ^L (a)}{b} \oplus \res{\Sigma ^R(a)}{b} =\\
& = \Sigma ^L(a) \oplus \Sigma ^R(b) = \res{\Sigma (b)}{a}
\end{split}
\end{equation*}
{\bf 7.ii.} \quad For a flag $ \lbrace b, x_j \rbrace$ of type $\{k+1, j \}$, with $j \in \{k+2, \dots,
k+4 \}$ we define:
\begin{center}
$\Sigma(b, x_j) := \res{\Sigma ^L(x_j)}{b} \oplus \Sigma ^R(x_j)$\\
\end{center}
Check the property $(S)$:
\begin{center}
$\res{\Sigma(b)}{x_j} = \res{\Sigma ^L (x_j)}{b} \oplus \Sigma ^R (x_j) = \res{\Sigma ^L (x_j) \oplus \Sigma ^R(x_j)}{b} = \res{\Sigma (x_j)}{b}$
\end{center}

Therefore $\Sigma$ is defined for all rank $1$ and rank $2$ flags of
$\Gamma$. Now using Ellard and Shult Lemma, we can extend the sheaf $\Sigma$
to a sheaf over all nonempty flags of $\Gamma$. We shall denote this sheaf
$\Sigma$, too, since this will not create confusion later.

\section{The proof of Theorem $3$}
In this section we shall prove the following:
\begin{tc} {\bf a}. Let {\pl} be a parapolar space which is locally $A_{n-1,3}({\mathbb K})$ where $n$ is an integer greater than $6$ and ${\mathbb K}$ is a field. Then $\Gamma$ is a residually connected locally truncated diagram geometry belonging to the diagram:\\
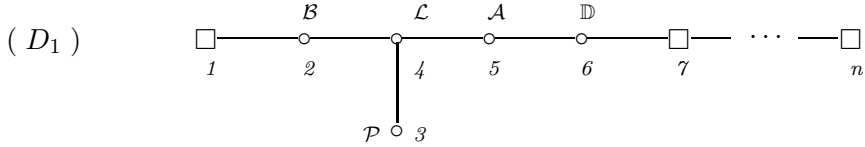
\begin{picture}(1000, 50)(0,0)
\put(0,25){$( \; D_1 \; )$}
\put(60,25){$\qed$}
\put(75,15){\scriptsize 1}
\put(80,29){\line(1,0){30}}
\put(110,26){$\circ$}
\put(112,36){\scriptsize ${\mathcal B}$}
\put(112,15){\scriptsize 2}
\put(115,29){\line(1,0){30}}
\put(145,26){$\circ$}
\put(154,36){\scriptsize ${\mathcal L}$}
\put(154,15){\scriptsize 4}
\put(148,-3){\line(0,1){30}}
\put(145,-9){$\circ$}
\put(135,-10){\scriptsize ${\mathcal P}$}
\put(154,-10){\scriptsize 3}
\put(150,29){\line(1,0){30}}
\put(180,26){$\circ$}
\put(182,36){\scriptsize ${\mathcal A}$}
\put(182,15){\scriptsize 5}
\put(185,29){\line(1,0){30}}
\put(215,26){$\circ$}
\put(217,36){\scriptsize ${\mathbb D}$}
\put(217,15){\scriptsize 6}
\put(220,29){\line(1,0){30}}
\put(239,25){$\qed$}
\put(253,15){\scriptsize 7}
\put(259,29){\line(1,0){15} \; $\dots$  \line(1,0){15}}
\put(304,25){$\qed$}
\put(319,15){\scriptsize n}
\end{picture}
\vspace{.3cm}

whose universal $2$-cover is the truncation of a building.\\
{\bf b.} Let {\pl} be a parapolar space which is locally $A_{n-1, 4}({\mathbb K})$ where $n=7$ or $n$ is an integer at least $9$ and ${\mathbb K}$ is a field. Assume that $\Gamma$ satisfies the Weak Hexagon Axiom. Then $\Gamma$ is a residually connected diagram geometry belonging to the diagram:\\
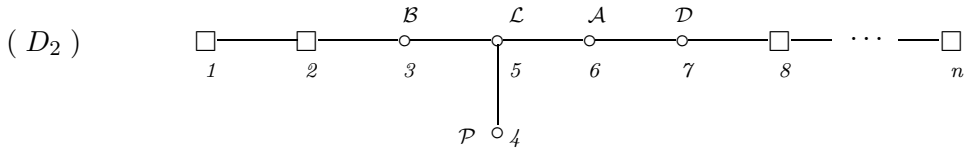
\begin{picture}(1000, 50)(0,0)
\put(0,25){$( \; D_2 \; )$}
\put(60,25){$\qed$}
\put(75,15){\scriptsize 1}
\put(80,29){\line(1,0){30}}
\put(98,25){$\qed$}
\put(113,15){\scriptsize 2}
\put(118,29){\line(1,0){30}}
\put(148,26){$\circ$}
\put(150,36){\scriptsize ${\mathcal B}$}
\put(150,15){\scriptsize 3}
\put(153,29){\line(1,0){30}}
\put(183,26){$\circ$}
\put(190,36){\scriptsize ${\mathcal L}$}
\put(190,15){\scriptsize 5}
\put(186,-3){\line(0,1){30}}
\put(183,-9){$\circ$}
\put(171,-10){\scriptsize ${\mathcal P}$}
\put(190,-10){\scriptsize 4}
\put(188,29){\line(1,0){30}}
\put(218,26){$\circ$}
\put(220,36){\scriptsize ${\mathcal A}$}
\put(220,15){\scriptsize 6}
\put(223,29){\line(1,0){30}}
\put(253,26){$\circ$}
\put(253,36){\scriptsize ${\mathcal D}$}
\put(255,15){\scriptsize 7}
\put(258,29){\line(1,0){30}}
\put(277,25){$\qed$}
\put(292,15){\scriptsize 8}
\put(297,29){\line(1,0){15} \; $\dots $ \line(1,0){15}}
\put(342,25){$\qed$}
\put(357,15){\scriptsize n}
\end{picture}
\vspace{.3cm}

whose universal $2$-cover is the truncation of a building.
\end{tc}
\begin{proof} Let $I = \lbrace 1, \ldots n \rbrace$ be a finite index set. Let $K \subset I$ and set $J = I \setminus K$. Assume {\pl} is a parapolar space, locally $A_{n-1,3}({\mathbb K})$ for $n$ an integer greater than $6$ and ${\mathbb K}$ a field. From the local properties, $\Gamma$ inherits two classes of maximal singular subspaces denoted ${\mathcal A}$ and ${\mathcal B}$. According to Theorem $1$, Section $4$, there exists a class of $2$-convex subspaces ${\mathbb D}$, each isomorphic to the half-spin geometry $D_{5,5}({\mathbb K})$. Therefore the parapolar space $\Gamma$ is enriched to a rank five geometry $\Gamma = ( {\mathcal B}, {\mathcal P}, {\mathcal L}, {\mathcal A}, {\mathbb D} )$ over $K = \lbrace 2, \ldots 6 \rbrace$. Then $\Gamma$ is a $J$-locally truncated geometry over $K$ belonging to diagram $(D_1)$ over $I$.

\medskip
Let now {\pl} satisfy the conditions of part $(b.)$ of the Theorem $3$. Therefore $\Gamma$ is a parapolar space, locally $A_{n-1, 4}({\mathbb K})$, with $n=7$ or an integer $n \geq 9$ and ${\mathbb K}$ is a field. There are two classes of maximal singular subspaces : ${\mathcal A}$ and ${\mathcal B}$. Also $\Gamma$ is assumed to satisfy the Weak Hexagon Axiom. Then there exists a family ${\mathcal D}$ of $2$-convex subspaces, each isomorphic to a $D_{6,6}({\mathbb K})$; see Theorem $2$, Section $5$. In this case $\Gamma$ is a rank five geometry $\Gamma =  ({\mathcal B}, {\mathcal P}, {\mathcal L}, {\mathcal A}, {\mathcal D} )$ over $K = \lbrace 3, \ldots 7 \rbrace$. Then $\Gamma$ is a $J$-locally truncated geometry over $K$ belonging to diagram $(D_2)$ over $I$.

\medskip
In what follows $\Gamma$ will stand for either of the two geometries mentioned above.
According to the previous section, there exists a sheaf $\Sigma$ of type
$I$ defined over the set ${\mathcal F}$ of all nonempty flags of $\Gamma$. For each flag $F \in {\mathcal F}$, the sheaf $\Sigma$ assigns a geometry $\Sigma(F)$ over $I \setminus t(F)$ with the property:
\begin{center}
$ _J \Sigma (F) \; = \; \res{\Gamma }{F}$
\end{center}
Since the sheaf values at rank $1$ and $2$ flags of $\Gamma$ are
$J$-truncated buildings or direct products of $J$-truncated buildings, the sheaf is residually connected.

\medskip
There exists a canonically defined chamber system ${\mathcal C}(\Sigma)$ over $I$. As proved by Brouwer and Cohen (Lemma $1$, \cite{loc}), this chamber system is residually connected and belongs to the diagram $(D_i), \; i=1,2$, also. It satisfies the property: \begin{center}
$_J {\mathcal C}(\Sigma) \; \simeq \; {\mathcal C}(\Gamma)$
\end{center}
where the above is an isomorphism of chamber systems over $K$ and with the term on the right being the chamber system corresponding to the geometry $\Gamma$.

\medskip
 The ambient diagram $(D_i)$, $i=1,2$ is a Coxeter diagram, thus of type $M$, and the chamber system ${\mathcal C}(\Sigma)$ has all its rank $3$ residues covered by truncations of buildings. Therefore by Tits' Local Approach Theorem \cite{tits} (see also Subsection $2.3$) the universal $2$-cover of ${\mathcal C}(\Sigma)$ is the chamber system ${\mathcal B}$ of a truncated building belonging to diagram $(D_i)$, $i=1,2$. Let $h: {\mathcal B} \rightarrow {\mathcal C}(\Sigma)$ be the $2$-covering map.

\medskip
Next apply the functor ${\mathbf G}: {\mathcal
Chamb}_I \rightarrow {\mathcal Geom}_I$ followed by the truncation functor $_J {\mathbf Tr}: {\mathcal Geom}_I \rightarrow {\mathcal Geom}_K$, see Subsection $2.3$, and get the following commutative diagram:
\begin{center}
$\begin{CD}
{\mathcal {B}}     @>{\mathbf{G}}>> \tilde\Delta
     @>_J{\mathbf{Tr}}>>_J\tilde\Delta\\
  @VV{h}V                 @VV{h_*}V
     @VV{h_{*J}}V\\
{\mathcal C}(\Sigma)          @>{\mathbf{G}}>>\Delta
   @>_J{\mathbf{Tr}}>>_J \Delta
\end{CD}$
\end{center}
\medskip
Here $\tilde \Delta = {\mathbf G}({\mathcal B})$ is a $J$-truncated building geometry belonging to
the diagram $D_i,\; i=1,2$ and ${\mathbf G}(h)=h_*: \tilde \Delta \rightarrow \Delta$ is the functorially
defined morphism of geometries induced by $h$. The chamber system ${\mathcal C}(\Sigma)$ is
residually connected, thus $\Delta =
{\mathbf G}({\mathcal C}(\Sigma))$, the geometry functorially defined by ${\mathcal C}(\Sigma)$,
is also residually connected. According to the above mentioned Lemma $1$ \cite{loc}, we have:
\begin{center}
$_J{\mathbf Tr}\; {\mathbf G}({\mathcal C}(\Sigma)) \simeq \; _J \Delta \simeq \Gamma$
\end{center}
Furthermore the functorially induced map on truncations:
\begin{center}
$h_{*J}: \; _J  \tilde \Delta  \rightarrow \; \Gamma$
\end{center}
is an epimorphism of geometries. This is a consequence of residual connectedness and of the fact that $h$ is a $2$-covering map in ${\mathcal Chamb}_I$, but the interested reader can find the details in Lemma $16$ \cite{ks}. This ends the proof of the Theorem.
\end{proof}
\bigskip
{\bf Acknowledgements}. I would like to thank to my advisor, Prof. Ernest Shult,
for suggesting this problem and for his continuous help and patience which made the completion of this work possible.


\begin{thebibliography}{20}
\bibitem{apt}A.E. Brouwer, A.M. Cohen, Some remarks on Tits geometries: with an
appendix by J. Tits, {\it Nederl. Akad. Weterisch. Indag. Mat.}, {\bf 45}(4)(1983),
393-402.
\bibitem{loc}A. Brouwer, A. Cohen, Local recognition of Tits geometries of classical
type, {\it Geom. Dedicata}, {\bf 20} (1986), 181-199.
\bibitem{hnbk}F. Buekenhout (editor) {\it Handbook of
Incidence Geometry}, Elsevier, (1995).
\bibitem{coh}A. Cohen, On a theorem of Cooperstein, {\it Europ. J. Combin.}, {\bf
4} (1983), 107-126.
\bibitem{cc} A. Cohen, B. Cooperstein, A characterization of some geometries of Lie type, {\it Geom. Dedicata}, {\bf 15}(1983), 73-105.
\bibitem{als}M. El-Atrash, E. Shult, A characterization of certain families of strong parapolar spaces, {\it Geom. Dedicata} {\bf 71} (1998), 221-235.
\bibitem{els}C. Ellard, E. Shult, A characterization of polar Grassmann spaces,
Preprint KSU (1981).
\bibitem{hans}G. Hanssens, A characterization of point-line geometries for finite
buildings, {\it Geom. Dedicata} {\bf 25} (1988), 297-315.
\bibitem{ks2}A. Kasikova, E. Shult, Chamber systems which are not geometric, {\it Comm. Algebra} {\bf 24} (1996), 3471-3481.
\bibitem{ks}A. Kasikova, E. Shult, Point-line characterizations of Lie
geometries, {\it Adv. Geom.} {\bf 2} (2002), no. 2, 147-188.
\bibitem{thesis}S. Onofrei, Characterizations of two classes of locally truncated diagram geometries, PhD Thesis, Kansas State University (2003).
\bibitem{pas}A. Pasini, {\it Diagram Geometries}, Oxford University Press (1994).
\bibitem{ron}M.A. Ronan, Extending locally truncated buildings and chamber systems,
{\it Proc. London Math. Soc.} 3, {\bf 53} (1986), 385-406.
\bibitem{shm}E.E. Shult, Characterizations of spaces related to metasymplectic spaces, {\it Geom. Dedicata} {\bf 30} (1989), no. 2, 325-371.
\bibitem{rem1}E. Shult, A remark on Grassmann spaces and half-spin geometries,
{\it Europ. J. Combinatorics} {\bf 15} (1994), 47-52.
\bibitem{asp}E. Shult, Aspects of buildings, Trends in Mathematics, Birkh\"auser (1998) 177-188.
\bibitem{ln}E. Shult, {\it Points and Lines}, unpublished book.
\bibitem{tits} J. Tits, A local approach to buildings. In {\it The Geometric
Vein}, The Coxeter Festschrift, Springer, Berlin (1981) 519-547.
\end{thebibliography}
\end{document}